\documentstyle[euscript,epsf,amssymb,12pt,xr]{amsart}

\input xy
\xyoption{all}

\newcommand{\la}{\langle}
\newcommand{\ra}{\rangle}

\newtheorem{theorem}{\bf Theorem}[section]
\newtheorem{lemma}[theorem]{\bf Lemma}

\newtheorem{corollary}[theorem]{\bf Corollary}
\newtheorem{definition}[theorem]{\bf Definition}

\newcommand{\CC}{{\Bbb C}}
\newcommand{\CP}{{\Bbb CP}}
\newcommand{\DD}{{\Bbb D}}

\newcommand{\JJ}{{\Bbb J}}

\newcommand{\NN}{{\Bbb N}}

\newcommand{\QQ}{{\Bbb Q}}
\newcommand{\RR}{{\Bbb R}}

\newcommand{\ZZ}{{\Bbb Z}}

\newcommand{\glie}{{\frak g}}

\newcommand{\tlie}{{\frak t}}

\newcommand{\CS}{{\operatorname {CS}}}
\newcommand{\bCS}{{\mathbf {CS}}}

\newcommand{\Diff}{\operatorname{Diff}}

\newcommand{\edge}{\operatorname{edge}}

\newcommand{\End}{\operatorname{End}}

\newcommand{\Hull}{\operatorname{Hull}}
\newcommand{\Id}{\operatorname{Id}}

\newcommand{\inter}{\operatorname{int}}

\newcommand{\Lie}{\operatorname{Lie}}
\newcommand{\length}{\operatorname{length}}

\newcommand{\PD}{\operatorname{PD}}
\newcommand{\pt}{\operatorname{pt}}

\newcommand{\rk}{\operatorname{rk}}

\newcommand{\textbegin}{\operatorname{begin}}
\newcommand{\textend}{\operatorname{end}}

\newcommand{\vertex}{\operatorname{vertex}}

\newcommand{\weight}{\operatorname{weight}}

\renewcommand{\exp}{\operatorname{exp}}

\newcommand{\ov}{\overline}

\newcommand{\curly}{\cal}

\newcommand{\CCC}{{\curly C}}

\newcommand{\III}{{\curly I}}
\newcommand{\JJJ}{{\curly J}}
\newcommand{\LLL}{{\curly L}}
\newcommand{\MMM}{{\curly M}}
\newcommand{\NNN}{{\curly N}}
\newcommand{\OOO}{{\curly O}}
\newcommand{\QQQ}{{\curly Q}}
\newcommand{\PPP}{{\curly P}}
\newcommand{\RRR}{{\curly R}}
\newcommand{\SSS}{{\curly S}}
\newcommand{\TTT}{{\curly T}}
\newcommand{\UUU}{{\curly U}}
\newcommand{\VVV}{{\curly V}}

\newcommand{\XXX}{{\curly X}}
\newcommand{\ZZZ}{{\curly Z}}
\newcommand{\ggreat}{>\kern-.7ex>}
\newcommand{\lless}{<\kern-.7ex<}
\newcommand{\qu}{/\kern-.7ex/}
\newcommand{\exh}{\to\kern-1.8ex\to}

\newcommand{\VP}{{\curly V}\kern-0.9ex\PPP}

\newcommand{\imag}{{\mathbf i}}

\newcommand{\bu}{{\mathbf{u}}}
\newcommand{\bx}{{\mathbf{x}}}

\newcommand{\bK}{{\mathbf{K}}}

\newcommand{\ZII}{Z_{I\kern-.3ex I}}
\newcommand{\ZIII}{Z_{I\kern-.3ex I\kern-.3ex I}}

\date{November 9, 2006}

\author[I. Mundet i Riera]{I. Mundet i Riera}
\address{Departament d'\`Algebra i Geometria, Facultat de Matem\`atiques,
Universitat de Barcelona, Gran Via de les Corts Catalanes 585,
08007 Barcelona, Spain} \email{ignasi.mundet@@ub.edu}

\title{The biinvariant diagonal class for Hamiltonian torus actions}

\begin{document}

\maketitle

\begin{abstract}
Suppose that an algebraic torus $G$ acts algebraically on a
projective manifold $X$ with generically trivial stabilizers. Then
the Zariski closure of the set of pairs $\{(x,y)\in X\times X\mid
y=gx \text{ for some }g\in G\}$ defines a nonzero equivariant
cohomology class $[\Delta_G]\in H^*_{G\times G}(X\times X)$. We
give an analogue of this construction in the case where $X$ is a
compact symplectic manifold endowed with a hamiltonian action of a
torus, whose complexification plays the role of $G$. We also prove
that the Kirwan map sends the class $[\Delta_G]$ to the class of
the diagonal in each symplectic quotient. This allows to define a
canonical right inverse of the Kirwan map.
\end{abstract}

\tableofcontents

\section{Introduction}

\label{s:intro}

\subsection{}
\label{ss:cas-algebraic}

The purpose of this paper is to give a symplectic version of the
following construction in algebraic geometry. Let $X$ be a smooth
projective scheme over $\CC$ of complex dimension $n$ endowed with
an algebraic action of an algebraic group $G$ with generically
trivial isotropy groups. Consider the following subset of $X\times
X$:
$$\Delta_G=\{(x,y)\in X\times X\mid \text{ there is some $g\in G$
such that $y=gx$ }\}.$$ The set $\Delta_G$ is a constructible
subset of $X\times X$ because it is the image of the algebraic map
$f:X\times G\to X\times X$ defined as $f(x,g)=(x,gx)$. This
implies that the inclusion of $\Delta_G$ in its Zariski closure
$\ov{\Delta}_G\subset X\times X$ has dense image with respect to
the analytic topology (see Ex. 3.18 and 3.19 in Chapter II of
\cite{Ha}, or \S 4.4 in \cite{Hu}) and hence that the dimension of
$\ov{\Delta}_G$ is equal to that of $\Delta_G$. Since the isotropy
groups of the action are generically trivial, the dimension of
$\Delta_G$ is equal to $n+\dim G$ (see for example Ex. 3.22 in
Chapter II of \cite{Ha}), so via the cycle map and Poincar\'e
duality $\ov{\Delta}_G$ defines a nonzero cohomology class
$[\Delta_G]_0\in H^{2(n-\dim G)}(X\times X;\ZZ)$. The set
$\ov{\Delta}_G$ is invariant under the product action of $G\times
G$ on $X\times X$ (in contrast with the usual diagonal which is
only invariant under the diagonal action of $G$ on $X\times X$).
Using algebraic finite dimensional approximations of the
classifying space $BG$ together with a stabilization argument, one
can apply the previous reasoning to define a nonzero equivariant
cohomology class
$$[\Delta_G]\in H^{2(n-\dim G)}_{G\times G}(X\times X;\ZZ).$$
We call $[\Delta_G]$ the biinvariant diagonal class.

\subsection{}
In this paper we generalize the previous construction to
symplectic geometry when $G$ is the complexification of a compact
torus $T$. We also prove some properties of $[\Delta_{G}]$, which
we construct over the rationals and  not over the integers. Let
$(X,\omega)$ be a compact connected symplectic manifold of real
dimension $2n$, endowed with an effective Hamiltonian action of
$T$, which for the moment we take to be $S^1$. Denote by $\mu:X\to
(\imag\RR)^*$ the moment map and define the function $h:X\to\RR$
as $h=\la\mu,\imag\ra$. Fix an invariant Riemannian metric on $X$
of the form $\omega(\cdot,I\cdot)$, where $I$ is an invariant
almost complex structure on $X$. Let $\xi_t:X\to X$ be the
downward gradient flow of $h$, defined by the conditions that
$\xi_0$ is the identity and $\xi_t'=-\xi_t^*\nabla h$. Define
$$\Delta_{\CC^*}=\{(x,y)\in X\times X\mid \text{ there is some $t\in\RR$
and $\theta\in S^1$ such that $y=\theta\cdot\xi_t(x)$ }\}.$$ When
$[\omega/2\pi]\in H^2(X;\ZZ)$ and $I$ is integrable then $X$ is
projective by Kodaira's theorem, the action of $S^1$ extends to an
algebraic action of $\CC^*$, and for any $z\in\CC^*$ and $x\in X$
we have $z\cdot x=\theta\cdot\xi_t(x)$, where $\theta=z/|z|$ and
$t=\ln|z|$. Hence this definition of $\Delta_{\CC^*}$ generalizes
the one in \S\ref{ss:cas-algebraic} and consequently, denoting by
$\ov{\Delta}_{\CC^*}$ the closure of $\Delta_{\CC^*}$ in the
standard topology of $X\times X$, the complement
$\ov{\Delta}_{\CC^*}\setminus {\Delta}_{\CC^*}$ has a natural
structure of stratified space of real dimension
$\dim_{\RR}\Delta_{\CC^*}-2$.

Unlike in the algebraic case, in general there is no reason to
expect that $\ov{\Delta}_{\CC^*}\setminus \Delta_{\CC^*}$ has
smaller dimension than $\Delta_{\CC^*}$ in any sense which would
allow $\ov{\Delta}_{\CC^*}$ to define a homology class of real
dimension $2n+2$. However, using multivalued perturbations of the
gradient flow equation (see for example \S 5.2 in \cite{S}, or \S
\ref{ss:multivalued-perturbations} in this paper for the notion of
multivalued perturbation) we can define a nonzero rational
cohomology class
$$[\Delta_{\CC^*}]\in H^{2n-2}_{S^1\times
S^1}(X\times X),$$ which is morally the equivariant Poincar\'e
dual of the class represented by $\Delta_{\CC^*}$, and which
coincides in the algebraic case with the class defined in
\S\ref{ss:cas-algebraic}. (Here and in the rest of the paper we
omit the coefficients in (co)homology groups, which are always
assumed to be $\QQ$.)

The idea of considering multivalued perturbations for the gradient
line equation to achieve simultaneously equivariance and
transversality has been well known to experts for some time. For
example, a sketch of this technique is explained in Lemma 4.7 of
\cite{McDT}, where it is applied to the definition of the
cohomology classes represented by stable and unstable manifolds.
However, as far as we know a detailed exposition of this
construction applied to the gradient line equation does not exist
in the literature, and this was one of the motivations for writing
this paper. Note that, in contrast, full details have been given
of the technique of multivalued perturbations applied to much more
involved geometric problems, such as the construction of
Gromov--Witten invariants or Floer homology \cite{FO,LiTi,R,S}.

\subsection{}
The main result of this paper is, besides the definition of
$[\Delta_{\CC^*}]$, the computation of its image under the
diagonal Kirwan map. Recall that if $m$ is a regular value of $h$
the symplectic quotient (or reduced space, or Marsden--Weinstein
quotient) of $X$ at $m$ is
$$Y_m=h^{-1}(m)/S^1.$$
The Kirwan map is the morphism of rings $$\kappa_m:H^*_{S^1}(X)\to
H^*(Y_m)$$ defined as the composition of the restriction map
$H^*_{S^1}(X)\to H^*_{S^1}(h^{-1}(m))$ with the Cartan isomorphism
$H^*_{S^1}(h^{-1}(m))\simeq H^*(Y_m)$, which exists because the
action of $S^1$ on $h^{-1}(m)$ has finite stabilizers. Now
$(X\times X)_{S^1\times S^1}$ can be identified with
$X_{S^1}\times X_{S^1}$, so K\"unneth's formula gives an
isomorphism
\begin{equation}
\label{eq:Kunneth} H^*_{S^1\times S^1}(X\times X)\simeq\bigoplus
H^*_{S^1}(X)\otimes H^*_{S^1}(X).
\end{equation}
Denote by $\kappa^2_m:H^*_{S^1\times S^1}(X\times X)\to
H^*(Y_m\times Y_m)$ the Kirwan map for the quotient of $X\times X$
at $(m,m)\in \RR^2$. Let $[\Delta_m]\in H^*(Y_m\times Y_m)$ denote
the Poincar\'e dual of the diagonal class (Poincar\'e duality
holds on $Y_m$ with rational coefficients because $Y_m$ is an
orbifold). We then have:

\begin{theorem}
\label{thm:classe-diagonal} For each regular value $m\in\RR$ of
$h$ we have $\kappa^2_m([\Delta_{\CC^*}])=[\Delta_m]$.
\end{theorem}

\subsection{A right inverse of the Kirwan map at regular quotients}
Let $m\in\RR$ be a regular value of the moment map. We deduce a
number of consequences of Theorem \ref{thm:classe-diagonal} by
looking at $(\kappa_m\otimes\Id) [\Delta_{\CC^*}]$ as a
correspondence in the sense of intersection theory. More
precisely, using Poincar\'e duality $PD:H^k(Y_m)\to
H^{2n-2-k}(Y_m)^*$, we can write
$$(PD\otimes\Id)\circ (\kappa_m\otimes\Id)[\Delta_{\CC^*}]\in
\bigoplus_{p+q=2n-2} H^{2n-2-p}(Y_m)^*\otimes H_{S^1}^q(X)
=\bigoplus_{q} H^{q}(Y_m)^*\otimes H_{S^1}^q(X).$$ Hence
$(PD\otimes\Id)\circ(\kappa_m\otimes\Id)[\Delta_{\CC^*}]$ gives
rise to a degree preserving linear map
$$l_m:H^*(Y_m)\to H^*_{S^1}(X).$$
Using Theorem \ref{thm:classe-diagonal} we prove the following.

\begin{corollary}
\label{cor:existence-inverse} The map $l_m$ is a left inverse of
the Kirwan map, i.e., the composition $\kappa_m\circ l_m$ is the
identity on $H^*(Y_m)$. In particular the Kirwan map $\kappa_m$ is
surjective.
\end{corollary}

The map $l_m$ is not in general morphisms of rings (see \S
\ref{ss:example} for an example).

To state the next corollary we need to introduce some notation.
The Poincar\'e dual of the class of the diagonal $[\Delta_m]\in
H^*(Y_m\times Y_m)$ gives rise to a nondegenerate quadratic $Q_m$
form on the homology $H_*(Y_m)$ defined as $Q_m(a,b)=\la a\otimes
b,[\Delta_m]\ra$ for any $a,b\in H_*(Y_m)$ (this is the usual
intersection product in homology). We can define similarly a
quadratic form $Q_{\CC}$ on the equivariant homology
$H^{S^1}_*(X)$ by setting $Q_{\CC}(\alpha,\beta)=\la
\alpha\otimes\beta,[\Delta_{\CC^*}]\ra$ for any $\alpha,\beta\in
H^{S^1}_*(X)$. Theorem \ref{thm:classe-diagonal} implies the
following.

\begin{corollary}
Let $m$ be any regular value of the moment map and let
$\kappa_m^*:H_*(Y_m)\to H_*^{S^1}(X)$ be the dual of the Kirwan
map. For any classes $a,b\in H_*(Y_m)$ we have
$$Q_{\CC}(\kappa_m^*(a),\kappa_m^*(b))=Q_m(a,b).$$
\end{corollary}

Note that the quadratic form $Q_{\CC}$ is always degenerate.

\subsection{Modified product in equivariant cohomology}
A way to encode the map $l_m$ is in terms of an associative ring
structure on the equivariant cohomology of $X$ which is different
from the usual one. Given classes $\alpha,\beta\in H^*_{S^1}(X)$
we define
$$\alpha\cup_m\beta=l_m(\kappa_m(\alpha)\cup\kappa_m(\beta)).$$
The product $\cup_m$ is associative because $l_m$ is a right
inverse of $\kappa_m$. In a sense, this is nothing but the usual
product in the cohomology of the symplectic quotient transported
via $l_m$ to the equivariant cohomology. What makes this
construction interesting is the possibility to define associative
deformations of $\cup_m$ in terms of the so-called Hamiltonian
Gromov--Witten invariants counting twisted holomorphic maps from
$\CP^1$ to $X$, similarly to how the quantum product is defined
(see \cite{MT1,MT2} and the references therein).

\subsection{The case of singular quotients}
\label{ss:singular-quotients} When $m$ is a critical value of $h$
the Kirwan map can not be defined as in the case of regular
values, since the cohomologies $H^*_{S^1}(h^{-1}(m))$ and
$H^*(Y_m)$ need no longer be isomorphic. In this situation, the
quotient $Y_m$ being a singular stratified space, it is more
natural to consider the (middle perversity) intersection
cohomology $IH^*(Y_m)$ rather than singular cohomology. Lerman and
Tolman have shown in \cite{LeTo} how to relate the equivariant
cohomology of $X$ to the intersection cohomology $IH^*(Y_m)$:
assuming that $m$ belongs to the interior of $h(X)$ (otherwise
$h^{-1}(m)$ is a connected component of the fixed point set), they
construct:
\begin{itemize}
\item an $S^1$-invariant Bott--Morse function $h':X\to\RR$, which
is a slight perturbation of $h$, such that $m$ is a regular value
of $h'$ and the action of $S^1$ on ${h'}^{-1}(m)$ has finite
stabilizers; \item a map $f:{h'}^{-1}(m)/S^1\to h^{-1}(m)/S^1$
which is a small resolution and hence induces an isomorphism
$f_H:IH^*(h^{-1}(m)/S^1)\simeq H^*({h'}^{-1}(m)/S^1)$ preserving
the intersection pairing.
\end{itemize}
Let $Y_m':={h'}^{-1}(m)/S^1$ and let us denote by
$\kappa_m':H^*_{S^1}(X)\to H^*(Y_m')$ the composition of the
restriction to ${h'}^{-1}(m)$ with the Cartan isomorphism
$H^*_{S^1}({h'}^{-1}(m))\to H^*(Y_m)$. It seems natural to call
the composition $\kappa_m=f_H^{-1}\circ \kappa_m'$ {\it the Kirwan
map} for the singular quotient $Y_m=h^{-1}(m)/S^1$. Kiem and Woolf
give in \cite{KW} a definition of Kirwan maps at singular
quotients for Hamiltonian actions of compact connected Lie groups.
Presumably the map $\kappa_m$ defined above can be obtained using
their technique (but note that the Kirwan maps constructed in
\cite{KW} are not canonical in general, whereas $\kappa_m$ is
canonical). Let us denote by $PD:IH^k(Y_m)\to IH^{2n-2-k}(Y_m)^*$
the Poincar\'e duality map.

\begin{theorem}
\label{thm:singular-m} The element $(PD\otimes\Id)\circ
(\kappa_m\otimes\Id)[\Delta_{\CC^*}]\in \bigoplus_{q}
IH^{q}(Y_m)^*\otimes H_{S^1}^q(X)$ corresponds a degree preserving
map $l_m:IH^*(Y_m)\to H_{S^1}^*(X)$ which is a right inverse of
the Kirwan map, i.e., $\kappa_m\circ l_m$ is the identity in
$IH^*(Y_m)$.
\end{theorem}

\subsection{Actions of any compact torus}
Now suppose that $X$ supports a Hamiltonian action of a compact
torus $T$ with Lie algebra $\tlie$ and moment map
$\mu:X\to\tlie^*$. Take a basis $u_1,\dots,u_q$ of $\tlie$ and
consider continuous curves in $X$ which are piecewise gradient
lines for $\la\mu,u_j\ra$ (see \S\ref{s:high-dim} for details).
Considering multivalued perturbations of the gradient line
equations which are invariant under the action of a generic
subgroup $S^1\simeq T_0\subset T$, one obtains a well defined
cohomology class, which is independent of the basis $\{u_i\}$ (but
maybe depends on the choice of $T_0$). One can then prove the
following.

\begin{theorem}
\label{thm:classe-diagonal-2} Let $q$ be the dimension of $T$. Let
$S^1\simeq T_0\subset T$ be a subgroup such that the $T_0$-fixed
point set coincides with the $T$-fixed point set. There is a
cohomology class
$$[\Delta_{\CC^*}^{T_0,T}]\in H^{2n-2q}_{T\times T}(X\times X)$$
such that, for each regular value $m\in\tlie^*$,
$\kappa_m^2([\Delta_{\CC^*}^{T_0,T}])=[\Delta_m].$
\end{theorem}

Similarly as in the case of $S^1$, the cohomology classes
$[\Delta_\CC^{T_0,T}]$ give rise to right inverses of the Kirwan
map
$$l_m^{T_0,T}:H^*(Y_m)\to H^*_T(X).$$

\subsection{Some remarks and questions}

If $T=S^1$ acts on $X$ is quasi-freely (i.e., all isotropy groups
are either trivial or the whole circle) then both
$[\Delta_{\CC^*}]$ and the right inverse $l_m$ can be defined over
the integers (using that the symplectic quotients are smooth and
hence that Poincar\'e duality holds over the integers). Moreover,
one can perturb the gradient flow equation using standard
perturbations (i.e., not multivalued), hence everything is much
easier. An immediate corollary is that if $H^*_T(X;\ZZ)$ is
torsion free then the cohomology of all symplectic quotients is
torsion free. This is a particular case of Theorem 5 in \cite{TW},
since for quasi-free actions of $S^1$ the group $H^*_T(F;\ZZ)$ is
torsion free if and only if $H^*_T(X;\ZZ)$ is torsion free (this
can be proved using the fact that, the action being quasi-free,
the moment map is a perfect Bott--Morse function over any finite
field).

An obvious question is whether the results in this paper extend to
other situations in which the Kirwan map is known to be
surjective: notably, the case of compact nonabelian groups
(already considered by Kirwan) and the case of loop group actions,
studied recently by Bott, Tolman and Weitsman in \cite{BTW}. More
generally, given any Hamiltonian action of a group $G$ on a
symplectic manifold $X$, one would like to understand the set of
cohomology classes $[\Delta_{G^{\CC}}]\in H^*_{G\times G}(X\times
X)$ such that for each regular value $\alpha\in\glie^*$ of the
moment map $\mu$, denoting by $Y_\alpha=\mu^{-1}(\alpha)/G$ the
symplectic reduction, one has
$\kappa^2_\alpha([\Delta_{G^{\CC}}])=[\Delta_\alpha]$, where
$[\Delta_\alpha]\in H^*(Y_\alpha\times Y_\alpha)$ is the diagonal
class. A particular case of this question is whether the classes
$[\Delta_{\CC^*}^{T_0,T}]$ constructed in Theorem
\ref{thm:classe-diagonal-2} depend on the choice of $T_0$. Another
question is whether the class $[\Delta_G]$ can be defined over the
integers, as is the case in the algebraic situation described in
\S \ref{ss:cas-algebraic}.

\subsection{Contents of the paper}
We now describe the contents of the remaining sections. In \S
\ref{s:perturbed-gd} we define the perturbed gradient segments
which will be used to define the biinvariant diagonal. The
biinvariant diagonal is defined, modulo Theorems \ref{thm:atlas}
and \ref{thm:definition-pairing}, in \S \ref{s:def-D-C}. The
proofs of Theorems \ref{thm:atlas} and
\ref{thm:definition-pairing} are given in \S
\ref{s:parametrizing}. In \S \ref{s:proof-thm-S1} we prove Theorem
\ref{thm:classe-diagonal}, Corollary \ref{cor:existence-inverse}
and Theorem \ref{thm:singular-m}. Finally, in \S \ref{s:high-dim}
we consider the case of higher dimensional compact tori and sketch
the proof of Theorem \ref{thm:classe-diagonal-2}.

\subsection{Acknowledgements} I wish to thank the referee for the
suggestion to consider the case of singular quotients in \S
\ref{ss:singular-quotients}.

\section{Perturbed gradient segments}
\label{s:perturbed-gd}

Recall that $I$ denotes an $S^1$-invariant almost complex
structure on $X$ which is compatible with $\omega$, and that we
consider on $X$ the Riemannian metric $\omega(\cdot,I\cdot)$. Let
$\XXX$ be the vector field on $X$ generated by the infinitesimal
action of $\imag\in\imag\RR\simeq\Lie S^1$. The function $h$
satisfies $dh=\iota_{\XXX}\omega$, so its gradient is $\nabla
h=I\XXX$. In order to define the cohomology class
$[\Delta_{\CC^*}]$ we consider generic $S^1$-invariant
perturbations of the downward gradient equation
$\gamma'=-I\XXX(\gamma)$. The possible presence of finite isotropy
groups forces us to consider multivalued perturbations in order to
preserve $S^1$-invariance, which is crucial to bound the dimension
of $\ov{\Delta}_{\CC^*}\setminus\Delta_{\CC^*}$ (see Lemma
\ref{lemma:clau}).

\subsection{$\epsilon$-perturbed gradient segments and some lemmata}

Let $c_1<\dots<c_r\in\RR$ be the critical values of $h$. Since the
moment map is locally constant on the fixed point set $F\subset X$
and $F$ coincides with the set of critical points of $h$, we have
$h(F)=\{c_1,\dots,c_m\}$. Choose real numbers $a_1,\dots,a_{m-1}$
satisfying
$$c_1<a_1<c_2<\dots<c_{m-1}<a_{m-1}<c_m$$
and define $Z\subset X$ to be the union of the submanifolds
$h^{-1}(a_1),\dots,h^{-1}(a_{m-1})$. Take a number $\beta>0$
satisfying $c_i+\beta<a_i<c_{i+1}-\beta$ for every $i$. Define
$$Z_i=h^{-1}([a_i-\beta,a_i+\beta])$$ for each $i$ and let $Z'$ be
the union $Z_1\cup\dots\cup Z_m$. Then the intersection $F\cap Z'$
is empty.

\begin{definition}
\label{def:perturbed-grad-segment} Let $A\subset\RR$ be an
interval and let $\epsilon$ be a positive number. A smooth map
$\gamma:A\to X$ is called an {\bf $\epsilon$-perturbed gradient
segment} if:
\begin{enumerate}
\item $\gamma'(t)=-I\XXX(\gamma(t))$ whenever $\gamma(t)\notin Z'$
and \item there exists a tangent vector field $\VVV$ defined in an
open neighborhood of the closure of $\gamma(A)\subset X$
satisfying $\gamma'=\VVV_{\gamma}$, $|\VVV+I\XXX|_{C^0}<\epsilon$
and $|\nabla(\VVV-I\XXX)|_{C^0}<\epsilon$.
\end{enumerate}
\end{definition}

Unless otherwise specified, the domain of an $\epsilon$-perturbed
gradient segment will always be assumed to be an interval of
$\RR$. Define the following quantities:
$$M=\sup_{Z'}|I\XXX|\qquad\text{ and }\qquad m=\inf_{Z'}|I\XXX|.$$
We will also always assume that $\epsilon\leq m/2$. This implies
that if $\gamma:A\to X$ is an $\epsilon$-perturbed gradient
segment then $(h\circ \gamma)'<0$, so $\gamma$ is and embedding
and the closure of $\gamma(A)$ is an embedded simply connected
curve.

\begin{lemma}
\label{lemma:length-osc} Suppose that $\epsilon<m^2$. Let
$\gamma:A\to X$ be an $\epsilon$-perturbed gradient segment and
assume that $B=\gamma^{-1}(Z_i)$ is nonempty. Then $B$ is
connected  and
\begin{equation}
\label{eq:length-osc} \frac{\length\gamma(B)}{M+\epsilon} \leq
|B|\leq \frac{\sup h(\gamma(B))-\inf h(\gamma(B))}{m(m-\epsilon)}.
\end{equation}
\end{lemma}
\begin{pf}
We first estimate for any $\gamma(t)\in Z'$, using $\nabla
h=I\XXX$ and Cauchy--Schwartz,
$$(h\circ\gamma)'(t)=\la\gamma'(t),\nabla h\ra=
-\la I\XXX,I\XXX\ra+\la\gamma'(t)+I\XXX,I\XXX\ra\leq -m^2+\epsilon
m.$$ We prove that $B$ is connected by contradiction. Suppose that
$t_0,t_1\in B$ but that $t_0<\tau<t_1$ satisfies
$\gamma(\tau)\notin Z_i$. Then either $h(\tau)>a_i+\beta$ or
$h(\tau)<a_i-\beta$. In the first case there must exist some
$t\in[\tau,t_1]$ such that $(h\circ\gamma)(t)=a_i+\beta$ and
$(h\circ\gamma)'(t)\geq 0$, which by the estimate above
contradicts our assumption $\epsilon<m$; the second case leads to
a contradiction in the same way. Once we know that $B$ is
connected, integrating the inequality $(h\circ\gamma)'(t)\leq
-m^2+\epsilon m$ along $B$ we get the second inequality in
(\ref{eq:length-osc}). To get the first inequality in
(\ref{eq:length-osc}) we estimate for any $\gamma(t)\in Z'$,
similarly as before, $|\gamma'(t)|\leq M+\epsilon$, and then we
integrate along $B$.
\end{pf}

For any $z\in Z$ and any real number $\delta>0$ we define the
following set:
$$S(z,\delta)=\{\exp_z v\mid v\in T_zX,
\ v\text{ is perpendicular to $\XXX$ and } |v|\leq\delta\}.$$ Let
also $R(z,\delta)=S^1\cdot S(z,\delta)$. If $\delta$ is smaller
than the injectivity radius then $S(z,\delta)$ is a slice of the
$S^1$ action at $z$. Also, if $\delta$ is small enough then
$R(z,\delta)$ can be identified with a solid torus (i.e., the
product of a circle and a closed ball) contained in $Z'$. The
generalized Gauss lemma for submanifolds (see e.g. Lemma 2.11 in
\cite{G}) implies that, given $x\in X$ and $z\in Z$,
\begin{equation}
\label{eq:Gauss} \text{$d(x,S^1\cdot z)$ is small enough
}\Longrightarrow d(x,S^1\cdot z)=\inf\{\delta\mid x\in
R(z,\delta)\}.
\end{equation}

\begin{lemma}
\label{lemma:convexitat} There exist positive numbers
$\epsilon,\delta_0$ with the following property. Suppose that
$\gamma:A\to X$ is an $\epsilon$-perturbed gradient segment. Let
$z\in Z$ and define, for any $t\in A$, $f(t)=\inf\{\delta^2\mid
\gamma(t)\in R(z,\delta)\}$. If $f(t)\leq \delta_0^2$, then
$f''(t)\geq 1$.
\end{lemma}
\begin{pf}
Take some $z\in Z$ and let $S\subset T_zX$ be the linear span of
$\XXX_z$ and $I\XXX_z$. Let $U\subset S^{\perp}$ be a small
neighborhood of $0$. Choose some small $\eta>0$ and let
$V=(-\eta,\eta)\times(-\eta,\eta)\times U$. The map $\iota:V\to X$
defined as $\iota(t,\theta,u)=\xi_t(e^{\imag\theta}\cdot\exp_zu)$
is an embedding and $\iota(V)$ is a neighborhood of $z$ in $X$
(recall that $\xi_t$ is the downward gradient flow of $h$). We can
identify $O=\{0\}\times(-\eta,\eta)\times\{0\}\subset V$ with
$\iota^{-1}(S^1\cdot z)$. Consider the Riemannian metric on $V$
defined by $dg^2=dt^2+d\theta^2+du^2$, where $du^2$ is the
restriction to $S^{\perp}$ of the Euclidean pairing in $T_zX$, and
let $d_0$ be the distance in $V$ induced by $dg^2$. The integral
curves of $d\iota^{-1}(I\XXX)$ on $V$ are given by
$\gamma(t):=(t,\theta_0,u_0)$ for some constant
$(\theta_0,u_0)\in(-\eta,\eta)\times U$, so the function
$f_0(t):=d_0(\gamma(t),O)^2$ satisfies $f''_0(t)=2$. Let $d$ be
the distance in $V$ induced by the distance in $X$ via the
inclusion $\iota$. If $\eta$ and $\epsilon$ are small enough and
$\VVV$ satisfies the hypothesis of Definition
\ref{def:perturbed-grad-segment}, then for any integral curve
$\gamma_{\VVV}$ of $\VVV$ the function
$f_{\VVV}=d(\gamma_{\VVV},O)^2$ satisfies $f_{\VVV}''\geq 1$. By
(\ref{eq:Gauss}) if $f$ is small enough then $f=f_{\VVV}$.
\end{pf}

\begin{lemma}
\label{lemma:e-petit} Let $\epsilon=m/2$. For any $\delta_0>0$
there is some $0<\delta_1<\delta_0$ such that if $\gamma:A\to X$
is an $\epsilon$-gradient segment and $z\in Z$ then
$\gamma(\Hull(\gamma^{-1}R(z,\delta_1)))\subset R(z,\delta_0)$,
where $\Hull$ denotes the convex hull.
\end{lemma}
\begin{pf}
Take a positive $\delta_2<\delta_0$ such that for any $z\in Z$ we
have $R(z,\delta_2)\subset Z'$. Let $d$ be the infimum for all
points $z\in Z$ of the distance between the boundaries $\partial
R(z,\delta_2/2)$ and $\partial R(z,\delta_2)$. If $\delta_2$ has
been chosen small enough, we have $d>0$. Choose a positive
$\delta_1<\delta_2/2$ in such a way that for any $z\in Z$ we have
\begin{equation}
\label{eq:torus-estret}
\sup_{R(z,\delta_1)}h-\inf_{R(z,\delta_1)}h < {dm^2}/{2M}.
\end{equation}
We prove that $\delta_1$ satisfies the requirement of the lemma,
even replacing $\delta_0$ by $\delta_2$. Suppose that $\gamma:A\to
X$ is an $\epsilon$-gradient segment, and that for some $z\in Z$
there exist elements $\tau<\tau'<\tau''$ of $A$ such that
$\gamma(\tau),\gamma(\tau'')\in R(z,\delta_1)$ but
$\gamma(\tau')\notin R(z,\delta_2)$. Let $B=[\tau,\tau'']$. By the
triangle inequality $\length(\gamma(B))\geq 2d$. Combining both
inequalities in (\ref{eq:length-osc}) we have
$h(\gamma(\tau))-h(\gamma(\tau''))\geq {dm^2}/{2M}$, contradicting
(\ref{eq:torus-estret}). This proves the lemma.
\end{pf}

\begin{lemma}
\label{lemma:e-petit-2} There exist positive numbers
$\epsilon,\delta$ with the following property. Suppose that
$\gamma:A\to X$ is an $\epsilon$-perturbed gradient segment. Let
$z\in Z$. Then
\begin{enumerate}
\item The preimage $\gamma^{-1}R(z,\delta)\subset A$ is connected.
\item If $\emptyset\neq\gamma(A)\cap R(z,\delta)\subset\partial
R(z,\delta)$ then $\gamma^{-1}R(z,\delta)$ consists of a unique
point.
\end{enumerate}
\end{lemma}
\begin{pf}
Let $\delta_0$ be given by Lemma \ref{lemma:convexitat} and let
$\delta_1$ be the corresponding value given by Lemma
\ref{lemma:e-petit}. Let $\epsilon$ be less than the $\epsilon$'s
in both lemmata and let $\delta=\delta_1$. If $\gamma:A\to X$ is
an $\epsilon$-perturbed gradient segment and $z\in Z$ then by
Lemma \ref{lemma:e-petit} the convex hull $B$ of
$\gamma^{-1}R(z,\delta)$ satisfies $\gamma(B)\subset
R(z,\delta_0)$. It follows that the function $f:B\to\RR$ defined
in Lemma \ref{lemma:convexitat} is convex, and hence
$\gamma^{-1}R(z,\delta)\cap B$ is connected. Hence,
$B=\gamma^{-1}R(z,\delta)$ and so the latter is connected. This
proves (1), and (2) follows similarly.
\end{pf}

\begin{lemma}
\label{lemma:tangencia} Suppose that $\gamma:A\to X$ is an
$\epsilon$-perturbed gradient segment and that $A$ is a closed
interval. Let $\VVV$ be the vector field defined in a neighborhood
of $\gamma(A)$ as given by Definition
\ref{def:perturbed-grad-segment}. Then the following is true. (1)
One can take open neighborhoods $H\subset\RR$ (resp. $O\subset X$)
of $h(\gamma(\sup A))$ (resp. $\gamma(\inf A)$) such that for any
$\lambda\in H$ and any $x\in O$ there is a unique integral curve
$\gamma_{\lambda,x}:A_{\lambda,x}\to X$ of $\VVV$ such that
$h(\gamma_{\lambda,x}(\sup A_{\lambda,x}))=\lambda$ and
$\gamma_{\lambda,x}(\inf A)=x$. (2) Take some point $z\in Z$ and
define, for small enough $\delta>0$, the following sets:
\begin{align*}
\Sigma_{\delta} &= \{(\lambda,x)\in H \times O\mid
\gamma_{\lambda,x}^{-1}R(z,\delta)\text{ consists of a unique point }\},\\
\Sigma_{\delta,1} &=\{(\lambda,x)\in \Sigma_{\delta} \mid
\gamma_{\lambda,x}^{-1}R(z,\delta)
\text{ belongs to the interior of $A_{\lambda,x}$}\}, \\
\Sigma_{\delta,2} &=\{(\lambda,x)\in \Sigma_{\delta} \mid
t=\gamma_{\lambda,x}^{-1}R(z,\delta)\in
\partial A_{\lambda,x} \text{ and $\gamma_{\lambda,x}'(t)$ is not tangent to
$\partial R(z,\delta)$}\}, \\
\Sigma_{\delta,3} &=\{(\lambda,x)\in \Sigma_{\delta} \mid
t=\gamma_{\lambda,x}^{-1}R(z,\delta)\in
\partial A_{\lambda,x} \text{ and $\gamma_{\lambda,x}'(t)$ is tangent to
$\partial R(z,\delta)$}\}
\end{align*}
Then $\Sigma_{\delta,3}\subset H\times O$ is a smooth submanifold
of dimension $2n-1$ and
$\Sigma_{\delta,1},\Sigma_{\delta,2}\subset H\times O$ are smooth
submanifolds of dimension $2n$. (2) If $j=1,2$ and
$p=(\lambda,x)\in\Sigma_{\delta,j}$, then there is an open
neighborhood $p\in\UUU\subset H\times O$ and an open interval
$D\subset\RR$ containing $\delta$ such that
$\{\UUU\cap\Sigma_{\delta',j}\}_{\delta'\in D}$ defines a smooth
foliation of $\UUU$. (3) Let $l:H\times O\to\RR$ be the map which
sends $(\lambda,x)$ to the length of
$\gamma_{\lambda,x}(A_{\lambda,x})\cap R(z,\delta)$. Then $l$ is a
continuous function.
\end{lemma}
\begin{pf}
Claim (1) follows from the existence and uniqueness of integral
curves of smooth vector fields. Statement (2) follows from
observing that $\VVV$ is tangent to $\partial R(z,\delta)$ along a
codimension $1$ submanifold of $\partial R(z,\delta)$. Finally,
(3) follows from the same arguments as in the proof of Lemma
\ref{lemma:convexitat}.
\end{pf}

\subsection{$J$-perturbed gradient segments}
\label{ss:multivalued-perturbations}

In this section we define multivalued perturbations of the
gradient flow equation $\gamma'=-I\XXX(\gamma)$ in terms of
infinitesimal variations of the almost complex structure. These
perturbations are called multivalued because they are defined on
finite nonramified coverings of the tori $R(z,\delta)$. Take a
point $z\in Z$ and suppose that the isotropy group of $z$ has $k$
elements. For any $\delta$ small enough so that $R(z,\delta)$ is a
solid torus and smaller than the injectivity radius we define
$$R^{\sharp}(z,\delta)=\{(\alpha,x)\in S^1\times R(z,\delta)\mid
\alpha^{-1}x\in S(z,\delta)\}.$$
Then the projection to the second factor
$$\pi:R^{\sharp}(z,\delta)\to R(z,\delta)$$
is an unramified covering of degree $k$ because the set of
elements $\theta\in S^1$ such that $\theta\cdot
S(z,\delta)=S(z,\delta)$ coincides with the stabilizer of $z$
(here we use that $\delta$ is less than the injectivity radius).
On the other hand, if we denote by $\OOO\subset R(z,\delta)$ the
$S^1$ orbit through $z$, the covering $\pi^{-1}(\OOO)\to\OOO$ is
isomorphic to the map $S^1\to S^1$ which sends $\theta$ to
$\theta^k$. It follows that $\pi^{-1}(\OOO)$ is connected, and
hence so is $R^{\sharp}(z,\delta)$. Consequently,
$R^{\sharp}(z,\delta)$ is a solid torus. Consider the action of
$S^1$ on $R^{\sharp}(z,\delta)$ defined as
$\theta\cdot(\alpha,x)=(\theta\alpha,\theta\cdot x)$ for any
$\theta\in S^1$. This action is free and, with respect to this
action, $\pi$ is equivariant.

Let $(V,\eta)$ be a symplectic vector space and let
$\JJJ\subset\End V$ be the set of complex structures $J\in\End V$
such that $\eta(\cdot,J\cdot)$ defines an Euclidean pairing. The
tangent space $T_J\JJJ$ can be identified with the space of
endomorphisms $j\in\End V$ satisfying $jJ+Jj=0$ and
$j+j^{*_{\eta}}=0$, where $j^{*_\eta}$ is the dual of $j$ with
respect to $\eta$. Hence the sections of the vector bundle
$$E=\{j\in \End TX \mid jI+Ij=0,\ j+j^{*_\omega}=0\}$$
can be identified with the infinitesimal deformations of $I$
as an almost complex structure compatible with $\omega$.
The following lemma is elementary and well known.

\begin{lemma}
\label{lemma:span-1} For any $0\neq v\in V$ and any $J\in\JJJ$ the
map $T_{J}\JJJ\ni j\mapsto jv\in V$ is onto.
\end{lemma}

Fix some point $z\in Z$. The previous lemma implies that we can
find $j_1,\dots,j_k\in E_z$ such that for any $v\in T_zX$ the
vectors $j_1(v),\dots,j_k(v)$ span $T_zX$. Choose $\delta_z>0$
smaller than the injectivity radius and the $\delta$ in Lemma
\ref{lemma:e-petit-2}, such that $R(z,\delta_z)$ is a solid torus
and such that there exist sections $J_1,\dots,J_k\in
C^{\infty}(S(z,\delta_z);E)$ satisfying: (1) $J_i(z)=j_i$ and (2)
for any $z'\in S(z,\delta_z)$ and tangent vector $v\in T_{z'}X$
the vectors $J_1(z')v,\dots,J_k(z')v$ span $T_{z'}X$. The pullback
vector bundle
$$\pi^*E\to R^{\sharp}(z,\delta_z)$$
admits a canonical lift of the $S^1$ action on $R^{\sharp}$. Since
such action is free and $\{1\}\times S(z,\delta)\subset
R^{\sharp}(z,\delta)$ is a slice, one can extend uniquely the
sections $J_1,\dots,J_k$ to equivariant sections
$J_1^{\sharp},\dots,J_k^{\sharp}$ of the vector bundle $\pi^*E$.
Denote by $$\JJ_z\subset
C^{\infty}(R^{\sharp}(z,\delta_z);\pi^*E)$$ the span of the
sections $J_1^{\sharp},\dots,J_k^{\sharp}$. Let
$\beta:\RR\to\RR_{\geq 0}$ be a smooth nonincreasing function
satisfying, for a small $\varepsilon>0$, $\beta(t)=1$ if
$t<\varepsilon$ and $\beta(t)=0$ if $t>1-\varepsilon$. For any
positive $\delta<\delta_z$ denote by
$\eta_{z,\delta}:R(z,\delta_z)\to \RR_{\geq 0}$ the unique
invariant function whose restriction to $S(z,\delta_z)$ satisfies
$\eta_{z,\delta}(\exp v)=\beta(|v|/\delta)$ for any $v\in T_zX$.

To state the following lemma we need to introduce some notation.
Let $\epsilon>0$ be a small real number, let $\tau>0$ and let
$\gamma:[0,\tau]\to X$ be an $\epsilon$-perturbed gradient segment
whose image intersects the interior of $R(z,\delta_z)$. Combining
Definition \ref{def:perturbed-grad-segment}, the deduction before
Lemma \ref{lemma:length-osc}, and Lemma \ref{lemma:e-petit-2}, we
deduce that there exists a connected and simply connected
neighborhood $U\subset R(z,\delta_z)$ of $\gamma(A)\cap
R(z,\delta_z)$ and a vector field $\VVV$ on $U$ which is tangent
to $\gamma$. Choose a lift $\sigma:U\to R^{\sharp}(z,\delta_z)$.
If $j\in\JJ_z$ is sufficiently near $0$, then one can define
$\gamma_j:[0,\tau]\to X$ by the properties $\gamma_j(0)=\gamma(0)$
and $\gamma_j'(t)= \VVV_{\gamma_j(t)}-\eta_{z,\delta}(j\circ
\sigma\circ \gamma_j)(t)\XXX_{\gamma_j(t)}$ (one only needs that,
for any $t\in [0,\tau]$, $\gamma_j(t)$ stays in $U$). Hence
$e(j)=\gamma_j(\tau)$ is defined for any $j$ contained in a small
neighborhood of $0$ in $\JJ_z$.

\begin{lemma}
\label{lemma:de-onto} If $\epsilon$ is small enough then the map
$e$ is differentiable and the differential $de(0):\JJ_z\to
T_{\gamma(\tau)}$ is onto.
\end{lemma}
\begin{pf}
It follows from standard results on ODE's and the definition of
$\JJ_z$, using for example the same coordinate charts as in the
proof of Lemma \ref{lemma:convexitat}.
\end{pf}

The union of the interiors of the sets $R(z,\delta_z)$ as $z$ runs
over all points in $Z$ contains $Z$, so by compactness one can
take points $z_1,\dots,z_s\in Z$ such that $Z$ is contained in the
union of the interiors of the sets $R(z_i,\delta_{z_i})$. Let
$R_i=R(z_i,\delta_{z_i})$, $\eta_i=\eta_{z_i,\delta_{z_i}}$ and
$R_{i}^{\sharp}=R^{\sharp}(z_i,\delta_{z_i})$ and
$\JJ_i=\JJ_{z_i}$. Assume that $\epsilon$ is small enough so that
Lemmata \ref{lemma:e-petit-2} and \ref{lemma:de-onto} hold true.
For any $j\in\JJ_i$ define
$\|j\|=\sup_{R_i^{\sharp}}|j|+\sup_{R_i^{\sharp}} |\nabla(\eta_i
j)|$ and let also
$$\JJ=\{(j_1,\dots,j_s)\in
\JJ_1\times\dots\times\JJ_r\mid \|j_i\|<\epsilon/s\text{ for each
$i$ } \}.$$ Let $J=(j_1,\dots,j_s)\in\JJ$. A $J$-{\bf perturbed
gradient segment} is a tuple
$(\gamma,\gamma_1^{\sharp},\dots,\gamma_s^{\sharp})$, where
$\gamma:A\to X$ is an $\epsilon$-perturbed gradient segment and
each $\gamma_i^{\sharp}:\gamma^{-1}R_i\to R_i^{\sharp}$ is a lift
of the restriction of $\gamma$, so that $\pi\circ
\gamma_i^{\sharp}=\gamma$ holds on $\gamma^{-1}R_i$, satisfying
the equation
\begin{equation}
\label{eq:pert-grad-segm} \gamma'=-(I+\sum_i
\eta_{i}j_i(\gamma_i^{\sharp}))\XXX_{\gamma}.
\end{equation}
If $\gamma^{-1}R_i$ is empty then $\gamma_i^{\sharp}$ and its
contribution in the differential equation can be ignored. Since
$J\in \JJ$ and the functions $\eta_i$ are everywhere $\leq 1$,
equation (\ref{eq:pert-grad-segm}) implies that $\gamma$ is an
$\epsilon$-perturbed gradient segment. Let $S^1$ act on
$J$-perturbed gradient segments as
$$\theta\cdot(\gamma,\gamma_1^{\sharp},\dots,\gamma_s^{\sharp})
=(\theta\cdot\gamma,\theta\cdot\gamma_1^{\sharp},\dots,
\theta\cdot\gamma_s^{\sharp}).$$ Note that
$\theta\cdot\gamma_i^{\sharp}$ is a lift of $\theta\cdot\gamma$
because the covering maps $\pi:R^{\sharp}_i\to R_i$ are $S^1$
equivariant, and equation (\ref{eq:pert-grad-segm}) is preserved
because $I$, $j_i$ and $\eta_i$ are $S^1$ invariant.

\subsection{}
\label{ss:def-oriented-chain} We define an {\bf oriented chain of
$J$-perturbed gradient segments} to be any tuple
$C=(K,K_1,\dots,K_s,b)$, where $K\subset X$ is a compact subset,
each $K_i\subset R_i^{\sharp}$ is a compact (possibly empty)
subset and $b\in K$, subject to the following conditions.
\begin{enumerate}
\item There is a continuous injective map $\rho:B\to X$, where
$B\subset\RR$ is a compact interval, which induces a homeomorphism
between $B$ and $K$. \item For each $i$, $B_i=\rho^{-1}R_i$ is
connected and $\pi:R_i^{\sharp}\to R_i$ induces a homeomorphism
between $K_i$ and $\rho(B_i)$ (the latter set is independent of
the parametrization $\rho$). \item The set $\rho^{-1}(F)$ is
finite and for each connected component $B'\subset
B\setminus\rho^{-1}(F)$ there is a $J$-perturbed gradient segment
$(\gamma:A\to X,\gamma_1^{\sharp},\dots,\gamma_s^{\sharp})$ and an
increasing homeomorphism $g:B'\to A$ such that $\gamma\circ
g=\rho|_{B'}$ and the image of $\gamma_i^{\sharp}$ coincides with
$K_i$. \item The point $b$ is either $\rho(\inf B)$ or $\rho(\sup
B)$ (this indicates the orientation of the chain of gradient
segments).
\end{enumerate} We call $b$ the {\bf beginning} of $C$. If
$b=\rho(\sup B)$ then we define the {\bf end} of $C$ to be
$e=\rho(\inf B)$. Otherwise we define $e=\rho(\sup B)$. Given a
chain $C=(K,K_1,\dots,K_s,b)$ we define $d(K_i)=\sup\{d(x,\partial
R_i^{\sharp})\mid x\in K_i\}$, where $d(x,\partial R_i^{\sharp})$
is defined using the pullback to $R_i^{\sharp}$ of the Riemannian
metric on $X$ (if $K_i$ is empty then we set $d(K_i)=0$). Define
the distance between two chains $C=(K,K_1,\dots,K_s,b)$ and
$C'=(K',K'_1,\dots,K'_s,b')$ as
\begin{equation}
\label{eq:distancia-pert-segm} d(C,C')=d_H(K,K')+d(b,b')+\sum
d_H(K_i,K'_i) d(K_i)d(K'_i),
\end{equation}
where $d_H$ denotes the Hausdorff distance between sets, in the
first summand using the Riemannian metric on $X$ and in the third
summand using its pullback to $R_i^{\sharp}$ via $\pi$.

\subsection{The space of oriented chains of perturbed gradient segments}
\label{ss:space} Denote by $\CCC_J$ the set of oriented chains of
$J$-perturbed gradient segments modulo the relation which
identifies two chains $C,C'$ whenever $d(C,C')=0$. Note that if
the chains $C=(K,K_1,\dots,K_s,b)$ and
$C'=(K',K'_1,\dots,K'_s,b')$ are different but $d(C,C')=0$, then
$K=K'$ and for any $i$ such that $K_i\neq K'_i$ both $K_i$ and
$K_i'$ are contained in the boundary $\partial R_i^{\sharp}$ and
hence, by Lemma \ref{lemma:e-petit-2}, consist of a unique point
each. Take on $\CCC_J$ the topology induced by the distance $d$
and define an action of $S^1$ on $\CCC_J$ componentwise:
$\theta\cdot(K,K_1,\dots,K_s,b)= (\theta\cdot K,\theta\cdot
K_1,\dots,\theta\cdot K_s,\theta\cdot b).$ One checks that this
action maps elements of $\CCC_J$ to elements of $\CCC_J$ using the
action of $S^1$ of $J$-perturbed gradient segments defined above.
The proof of the following lemma is straightforward.

\begin{lemma}
\label{lemma:T-compact} $\CCC_J$ is compact, the action of $S^1$
on $\CCC_J$ is continuous, and the map
$$(b,e):\CCC_J\to X\times X$$
given by sending each $C\in\CCC_J$ to its beginning and end is
continuous.
\end{lemma}

\section{Definition of the class $[\Delta_{\CC^*}]$}
\label{s:def-D-C}

\subsection{}
For any $J\in\JJ_{\epsilon}$ define $\CCC_J^0\subset\CCC_J$ as the
set of chains of perturbed gradient segments $(K,K_1,\dots,K_s,b)$
such that $K\cap F=\emptyset$. We say that $C=(K,K_1,\dots,K_s,b)$
is tangent to $R_i$ if $\emptyset\neq K\cap R_i\subset\partial
R_i$ which implies by Lemma \ref{lemma:e-petit-2} that $K\cap R_i$
is one point. Define $\CCC_J^{0,0}$ as the set of chains
$C\in\CCC_J^0$ which are not tangent to any $R_i$. Let $o_i$ be
the degree of the covering $\pi:R_i^{\sharp}\to R_i$
(equivalently, the order of the isotropy group of $z_i$). We
define
$$\weight:\CCC_J^{0,0}\to\QQ$$
by sending $C=(K,K_1,\dots,K_s,b)\in\CCC_J^{0,0}$ to the product
$o_{i_1}^{-1}\dots o_{i_\nu}^{-1}$, where $\{i_1,\dots,i_{\nu}\}$
is the set of $i$ such that $K\cap R_i\neq\emptyset$ (since
$C\in\CCC_J^{0,0}$ this implies that $K\cap R_i$ contains points
in the interior of $R_i$). The next theorem will be proved in \S
\ref{ss:proof-thm:atlas}.

\begin{theorem}
\label{thm:atlas} For any $C\in\CCC_J^0$ there exist oriented
connected manifolds $U_1,\dots,U_N$ of real dimension $2n+1$ and
continuous maps $\phi_j:U_j\to\CCC_J^0$ satisfying these
properties:
\begin{enumerate}
\item For any $j$ the map $\eta_j:U_j\to\RR\times X$ which sends
$u\in U_j$ to $(h\circ e\circ \phi_j(u),b\circ \phi_j(u))$ is a
local diffeomorphism preserving the orientation. \item The union
$\phi_1(U_1)\cup\dots\cup\phi_N(U_N)$ is a neighborhood of $C$ in
$\CCC_J^0$.  \item If $C\in\CCC_J^{0,0}$ then $N$ can be taken to
be $1$.
\end{enumerate}
\end{theorem}

\subsection{}
We briefly recall the notion of pseudocycle introduced in \S 6.5
of \cite{McDS}. Let $N$ be a smooth manifold. A subset $R\subset
N$ is said to have dimension at most $d$ if there is a
$d$-dimensional manifold $S$ and a smooth map $g:S\to N$ such that
$R\subset g(S)$. Given a smooth map $f:M\to N$ of oriented
manifolds, the omega limit set of $f$, denoted $\Omega_f\subset
N$, is the intersection of all closed subsets $\ov{f(M\setminus
K)}\subset N$ as $K$ runs over the collection of all compact
subsets of $M$. If $M$ has dimension $d$, the map $f:M\to N$ is
called a $d$-dimensional pseudocycle if $\Omega_f$ has dimension
at most $d-2$. Two $d$-dimensional pseudocycles $f:M\to N$ and
$f':M'\to N$ are called bordant if there is an oriented manifold
$W$ of dimension $d+1$ with boundary $\partial W=M\cup (-M')$ and
a smooth map $F:W\to N$ extending $f$ and $f'$ such that
$\Omega_F$ has dimension at most $d-1$. In Remark 6.5.3 of
\cite{McDS} a construction is given which assigns to any
$d$-dimensional homology class $\beta\in H_d(N)$ a bordism class
of $d$-dimensional pseudocycles $f:M\to N$. We say that the
$f:M\to N$ represents $\beta$.

Recall that two smooth maps $\alpha:M\to N$ and $\alpha':M'\to N$
are said to be transverse if the map $(\alpha,\alpha'):M\times
M'\to N\times N$ is transverse to the diagonal $\Delta_N\subset
N\times N$. In this situation, the set
$\CS(\alpha,\alpha'):=\{(x,x')\in M\times M'\mid
\alpha(x)=\alpha'(x')\}$ is a submanifold of $M\times M'$ of
dimension $\dim M+\dim M'-\dim N$ (here $\CS$ stands for Cartesian
Square). If $\alpha:A\to B$ is a submersion, then $\alpha$ is
transverse to any smooth map $\alpha':A'\to B$. The proof of the
following lemma is straightforward.

\begin{lemma}
\label{lemma:CS-compacte} Suppose that $\alpha:M\to N$ and
$\alpha':M'\to N$ are two transverse maps satisfying
$\Omega_{\alpha}\cap\alpha'(M')=\alpha(M)\cap\Omega_{\alpha'}=\emptyset$.
Then $\CS(\alpha,\alpha')$ is a compact submanifold of $M\times
M'$.
\end{lemma}

\subsection{}
\label{ss:P-space} Define the set $\PPP_J=S^1\times \CCC_J$ and
let
\begin{equation}
\label{eq:map-Phi} \Theta_J:\PPP_J\to X\times X
\end{equation}
be the map $\Theta_J(\theta,C)=(\theta\cdot b(C),e(C))$.
Considering the action of $S^1\times S^1$ on $\PPP_J$ given by
$(\alpha,\beta)\cdot(\theta,C)=(\alpha\beta^{-1}\theta,\beta\cdot
C)$ and the product action on $X\times X$, the map $\Theta_J$ is
$S^1\times S^1$ equivariant. For any natural number $\Lambda$ let
$S_{\Lambda}$ be the unit sphere in $\CC^{\Lambda+1}$ centered at
the origin. Scalar multiplication gives a free action of $S^1$ on
$S_{\Lambda}$ and hence a structure of principal circle bundle on
the quotient map $S_\Lambda\to\CP^\Lambda=S_{\Lambda}/S^1$. The
bundles $S_{\Lambda}\to\CP^{\Lambda}$ provide finite dimensional
approximations of the universal circle fibration. Define
$E_{\Lambda}=S_{\Lambda}\times S_{\Lambda}$ and
$B_{\Lambda}=\CP^{\Lambda}\times\CP^{\Lambda}$. The natural
projection $p:E_{\Lambda}\to B_{\Lambda}$ endows $E_{\Lambda}$
with a structure of principal $S^1\times S^1$ bundle. Since
$\Theta_J$ is equivariant, it induces a map
\begin{equation}
\label{eq:def-Theta-J}
\Theta_{J,\Lambda}:\PPP_{J,{\Lambda}}=E_{\Lambda}\times_{S^1\times
S^1}\PPP_J\to X^2_{\Lambda}:=E_{\Lambda}\times_{S^1\times
S^1}(X\times X).
\end{equation}
The manifolds $B_{\Lambda}$ and $X\times X$ have natural
orientations, which induce an orientation on $X^2_{\Lambda}$. We
define a cohomology class $[\Delta_{\CC^*}]_{\Lambda}\in
H^{2n-2}(X^2_{\Lambda})$ in terms of its pairing with homology
classes. The following theorem will be proved in \S
\ref{ss:proof-thm:definition-pairing}.

\begin{theorem}
\label{thm:definition-pairing} Let $\beta\in
H_{2n-2}(X^2_{\Lambda})$, and let $f:M\to X^2_{\Lambda}$ be a
pseudocycle representing $\beta$. Let $D\subset
C^{\infty}(X^2_{\Lambda},TX^2_{\Lambda})$ be a linear subspace
such that for any $p\in X^2_{\Lambda}$ the evaluation map $D\to
T_pX^2_{\Lambda}$ is onto, and let $\DD=\{\exp \gamma\mid
\gamma\in D\}\subset \Diff(X^2_{\Lambda})$. There exists a
residual subset $\RRR\subset \JJ\times\DD$ such that for any
$(J,\xi)\in\RRR$ we have:
\begin{enumerate}
\item The set $\TTT_{J,\xi}=\{(x,y)\in \PPP_{J,\Lambda}\times
M\mid \Theta_{J,\Lambda}(x)=\xi\circ f(y)\}$ is finite. \item For
any $(x,y)\in\TTT_{J,\xi}$ we have $x\in
E_{\Lambda}\times_{S^1\times S^1}(S^1\times\CCC_J^{0,0})$. \item
Let $(x,y)\in\TTT_{J,\xi}$ and let $b=p(x)$. Let $O\subset
B_{\Lambda}$ be a small neighborhood of $b$. Take a trivialization
of $E_{\Lambda}|_O$ and denote by $\psi:O\times S^1\times\CCC_J\to
\PPP_{J,\Lambda}|_O$ the induced homeomorphism. Suppose that
$x=\psi(b,\theta,C)$. Let $\phi:U\to\CCC_J$ be a continuous map as
given by Theorem \ref{thm:atlas}, where $U$ is an oriented
$2n+1$-dimensional manifold and $\phi(U)$ is a neighborhood of
$C$. Endow $V=O\times S^1\times U$ with its product orientation,
and let $\phi_O:V\to O\times S^1\times\CCC_J$ be the map
$(\Id_O,\Id_{S^1},\phi)$. The differential $\delta$ at
$((b,\theta,C),y)$ of the map
$$(\Theta_J\circ\psi\circ\phi_O,\xi\circ f):V\times M\to X^2_{\Lambda}$$
is an isomorphism of vector spaces. Define $\sigma(x,y)=1$ if
$\delta$ preserves the orientations and $\sigma(x,y)=-1$
otherwise. Define also $\weight(x)=\weight(C)$. \item The
following number only depends on $\beta$ and $\Lambda$, and not on
$D,f,J,\xi$:
$$\Delta_{\Lambda}(\beta)=\sum_{(x,y)\in\TTT_{J,\xi}}
\sigma(x,y)\weight(x)\in\QQ.$$
\end{enumerate}
\end{theorem}

\subsection{Definition of the biinvariant diagonal class}
\label{ss:def-Delta} The map
$\Delta_{\Lambda}:H_{2n-2}(X^2_{\Lambda})\to\QQ$ defined by the
previous theorem is clearly linear and hence is induced by a
cohomology class $[\Delta_{\CC^*}]_{\Lambda}\in
H^{2n-2}(X^2_{\Lambda})$. To compare this class for different
values of $\Lambda$, note that there is a natural homotopy class
of inclusion $B_{\Lambda}\subset B_{\Lambda+1}$ whose image is the
product of two hyperplanes. This inclusion induces
$\iota_{\Lambda}:X^2_{\Lambda}\to X^2_{\Lambda+1}$. The same ideas
as in the proof of (4) of Theorem \ref{thm:definition-pairing}
imply that
$$\iota_{\Lambda}^{2n-2}[\Delta_{\CC^*}]_{\Lambda+1}=
[\Delta_{\CC^*}]_{\Lambda}.$$ For big enough $\Lambda$ the map
$\iota_{\Lambda}^{2n-2}$ is an isomorphism and we can identify the
cohomology groups $H^{2n-2}(X^2_{\Lambda})\simeq
H^{2n-2}_{S^1\times S^1}(X\times X)$. Hence the class
$[\Delta_{\CC^*}]_{\Lambda}$ defines an equivariant cohomology
class
$$[\Delta_{\CC^*}]\in H^{2n-2}_{S^1\times S^1}(X\times X).$$
We call $[\Delta_{\CC^*}]$ the biinvariant diagonal class.

\section{Parametrizing oriented $J$-perturbed chains of gradient
segments} \label{s:parametrizing}

\subsection{Proof of Theorem \ref{thm:atlas}}
\label{ss:proof-thm:atlas}

Let $C=(K,K_1,\dots,K_s,b)\in\CCC_J^0$ and assume that the
$J$-perturbed gradient segment $(\gamma:A\to
X,\gamma_1^{\sharp},\dots,\gamma_s^{\sharp})$ parameterizes $C$,
so $A\subset\RR$ is a compact interval, $K=\gamma(A)$ and
$K_i=\gamma_i^{\sharp}(\gamma^{-1}R_i)$. Assume that
$b=\gamma(\sup A)$. Let $\III'$ be the set of $i$'s such that $K$
intersects the interior of $R_i$ and let $\III''$ be the set of
$i$'s such that $C$ is tangent to $R_i$. Choose for any
$i\in\III'$ a small open neighborhood $M_{i}^{\sharp}\subset
R_i^{\sharp}$ of $K_{i}$ such that $\pi:M_{i}^{\sharp}\to
M_{i}:=\pi(M_{i}^{\sharp})$ is a diffeomorphism of open manifolds
with boundary, and let $\sigma_{i}:M_{i}\to M_{i}^{\sharp}$ be its
inverse. For any $i\in\III''$ let $q_i=K\cap R_i$, which by Lemma
\ref{lemma:e-petit-2} consists of a unique point, and let
$Q_i=\pi^{-1}(q_i)\subset R_i^{\sharp}$, which consists of $o_i$
different points. Let $\QQQ=\prod_{i\in\III''}Q_i$. Given
$q=(q_i)\in\QQQ$, choose for each $i\in\III''$ a small open
neighborhood $M_{i}^{\sharp}\subset R_i^{\sharp}$ of $q_i$ such
that $\pi:M_i^{\sharp}\to M_i:=\pi(M_i^{\sharp})$ is a
diffeomorphism of open manifolds with boundary, and let
$\sigma_i:M_i\to M_i^{\sharp}$ be its inverse. Let
$\III=\III'\cup\III''$. For any $q\in\QQQ$ let $M_q\subset X$ be
an open neighborhood of $\gamma(A)$ such that $M_q\cap R_i\subset
M_i$ for any $i\in\III$. Define the following vector field on
$M_q$
\begin{equation}
\label{eq:grad-perturbat} \VVV_q=-(I+\sum_{i\in\III}\eta_i
j_i(\sigma_i))\XXX.
\end{equation}
Then $\gamma:A\to X$ is an integral curve of $\VVV_q$.
Furthermore, all the integral curves of $\VVV_q$ satisfy the
conditions of Definition \ref{def:perturbed-grad-segment}. Let
$\lambda_0=h(\gamma(\inf A))$ and $x_0=\gamma(\sup A)$. Applying
Lemma \ref{lemma:tangencia} to $\VVV_q$ we obtain an open
neighborhood $U_q=H_q\times O_q\subset\RR\times X$ of
$(\lambda_0,x_0)$ and for each $(\lambda,x)\in U_q$ an integral
curve $\gamma_{\lambda,x}:A_{\lambda,x}\to X$, which is an
$\epsilon$-perturbed gradient segment. Let
$K_{\lambda,x}=\gamma_{\lambda,x}(A_{\lambda,x})$ and let
$K_{\lambda,x,i}=K_{\lambda,x}\cap R_i$. Taking $U_q$ small enough
we can assume that $K_{\lambda,x,i}$ is nonempty if and only if
$i\in\III$. Define
\begin{equation}
\label{eq:def-phi-q}
\phi_q(\lambda,x)=(K_{\lambda,x},K_{\lambda,x,1},\dots,K_{\lambda,x,s},
\gamma_{\lambda,x}(\sup A_{\lambda,x})).
\end{equation}
Then $\phi_q(\lambda,x)\in\CCC_J^0$, because it can be
parametrized by the $J$-perturbed gradient segment
$(\gamma_{\lambda,x},\gamma_{\lambda,x,1},\dots,\gamma_{\lambda,x,s})$,
where $\gamma_i^{\sharp}:\gamma_{\lambda,x}^{-1}R_i\to
R_i^{\sharp}$ is equal to $\sigma_i\circ \gamma_i$. In this way we
have defined a continuous map $\phi_q:U_q\to\CCC_J^0$. Picking the
right orientation of $U_q$ claim (1) of Theorem \ref{thm:atlas}
holds trivially. We prove that $\bigcup_{q\in\QQQ}\phi_q(U_q)$ is
a neighborhood of $C$ in $\CCC_J^0$, which is claim (2) of the
theorem. By (1) in Lemma \ref{lemma:e-petit-2} for any
$C'=(K',K_1',\dots,K_s',b')\in\CCC_J^0$ each intersection $K'\cap
R_i$ is connected, so if $C'$ lies sufficiently near $C$ the
compact $K'$ must be an integral curve of one of the vector fields
$\VVV_q$. When $C\in\CCC_J^{0,0}$ the set $\III''$ is empty, so
there is a unique open set $U$ and map $\phi:U\to\CCC_J^{0,0}$
whose image is a neighborhood of $C$. This proves (3). Finally, to
deal with the case $b=\gamma(\inf A)$ we proceed exactly as
before, replacing the last entry in (\ref{eq:def-phi-q}) by $x$.

\subsection{Proof of Theorem \ref{thm:definition-pairing}}

\label{ss:proof-thm:definition-pairing}

Let $\CCC= \{(J,C)\mid J\in\JJ,\ C\in\CCC_J\}$. Define the
distance between points in $\CCC$ as
$d((J,C),(J',C'))=\|J-J'\|+d(C,C'),$ where if $J=(j_1,\dots,j_s)$
and $J'=(j'_1,\dots,j'_s)$ then
$\|J-J'\|=\|j_1-j'_1\|+\dots+\|j_s-j'_s\|$ and $d(C,C')$ is
defined as in (\ref{eq:distancia-pert-segm}). Consider on $\CCC$
the topology induced by this distance and define the maps
$$(b,e):\CCC\to X\times X$$
by mapping $(J,C)\in\CCC$ to $(b(C),e(C))$, where $b(C),e(C)$ are
defined in \S \ref{ss:def-oriented-chain}. Consider also the
projection
$$\pi_J:\CCC\to\JJ$$
sending any $(J,C)\in\CCC$ to $J$. For any integer $r$ let
$\CCC^r=\{(J,(K,\dots))\in\CCC\mid \sharp K\cap F=r\}$ be the set
of perturbed chains which meet the fixed point set at $r$ points.
Fix from now on an orientation of $\JJ$.

\begin{lemma}
\label{lemma:no-critical-value-parametrizing} Let
${\bK}=(J,(K,\dots))\in\CCC^0$ and let $\LLL=\{l\mid K\cap\inter
R_l=\emptyset,\ K\cap
\partial R_l\neq\emptyset\}.$
For any $l\in\LLL$ let $q_l\in R_l$ be the unique point of
intersection of $K$ with $R_l$ (see (2) in Lemma
\ref{lemma:e-petit-2}), and let $Q_l\subset R_l^{\sharp}$ be the
preimage of $q_l$. Let $\QQQ=\prod_{l\in\LLL}Q_l$.
\begin{enumerate}
\item There exists a collection of connected oriented open
manifolds $\{\UUU_q\}_{q\in\QQQ}$ of dimension equal to
$2n+1+\dim\JJ$ and continuous maps $\Phi_q:\UUU_q\to\CCC$ such
that the union $\bigcup_{q\in\QQQ}\Phi_q(\UUU_q)$ is a
neighborhood of $\bK$ in $\CCC$.

\item For any $q$ both $\pi_J\circ\Phi_q:\UUU_q\to\JJ$ and
$(b,e)\circ\Phi_q:\UUU_q\to X\times X$ are smooth maps.

\item For any $l\in\LLL$ define
$\OOO_l:=\{(J,(K,\dots))\in\CCC\mid K\cap\inter R_l=\emptyset\}$.
Define also, for any $q$, $\OOO_{q,l}=\Phi_q^{-1}(\OOO_l)$. If
$q=(q_l)\neq q'=(q_l')$, then we have
$\Phi_q^{-1}(\Phi_{q'}(\UUU_{q'}))=\bigcap_{q_l\neq q_l'}
\OOO_{q,l}$. The boundary $\partial \OOO_{q,l}\subset\UUU_q$ is
the disjoint union of smooth submanifolds
$\SSS_{q,l,1},\SSS_{q,l,2},\SSS_{q,l,3}$ of codimensions $1$, $1$
and $2$ respectively.
\end{enumerate}
\end{lemma}
\begin{pf}
(1) and (2) follow from the same arguments as the proof of Theorem
\ref{thm:atlas} given in \S \ref{ss:proof-thm:atlas}, replacing
$X$ by $\JJ\times X$ and choosing the perturbations $j_i$ in
(\ref{eq:grad-perturbat}) using the coordinate in $\JJ$ (note that
$\LLL$ corresponds to $\III''$). The first statement in (3)
follows from the construction; in the second statement, the
submanifolds $\SSS_{q,j,i}$ are the analogues of the submanifolds
$\Sigma_{\delta,i}$ in Lemma \ref{lemma:tangencia}.
\end{pf}

\begin{lemma}
\label{lemma:broken-lines} Let $r\geq 1$ be an integer. There is a
countable collection of connected smooth manifolds
$\{\VVV_{r,i}\}_{i\in\NN}$ of dimension $2n+1-r$ and continuous
maps $\Psi_{r,i}:\UUU_{r,i}\to\CCC^r$ such that: (1) the union
$\bigcup_{i\in\NN}\Psi_{r,i}(\VVV_{r,i})$ is equal to $\CCC^r$,
(2) for each $i$ the compositions
$(b,e)\circ\Psi_{r,i}:\VVV_{r,i}\to X\times X$ and
$\pi_J\circ\Psi_{r,i}:\VVV_{r,i}\to\JJ$ are smooth maps.
\end{lemma}
\begin{pf}
Given a closed interval $[u,v]\subset\RR$ we define
$\CCC([u,v])\subset\CCC$ as the subset of all $\bK\in\CCC$ such
that $h(b(\bK))=u$ and $h(e(\bK))=v$. More generally, for any
interval $A\subset\RR$, let $\CCC(A)$ to be the union of all the
sets $\CCC(B)$ as $B$ runs over the collection of the compact
subintervals of $A$. Define also for any $r$ the set
$\CCC^r(A)=\CCC^r\cap\CCC(A)$. We prove the lemma in several
steps.

\noindent{\bf Step 1.} Let $\CCC^{\geq r}=\bigsqcup_{r'\geq
r}\CCC^{r'}$. It is straightforward to check (as in Lemma
\ref{lemma:T-compact}) that the projection $\pi_J:\CCC^{\geq
r}\to\JJ$ is proper. Furthermore, $\CCC^r$ is open in $\CCC^{\geq
r}$ so, defining for any integer $\alpha$ the subset
$\CCC^{r,\alpha}=\{\bK\in\CCC^r\mid d(\bK,\CCC^{\geq r+1})\in
[2^{-\alpha},2^{-\alpha+1}]\}\subset\CCC^r$, the restriction of
$\pi_J$ to each $\CCC^{r,\alpha}$ is proper, and also
$\CCC^r=\bigcup_{\alpha}\CCC^{r,\alpha}$. Hence it suffices to
construct for any $\bK\in\CCC^r$ a collection of connected
manifolds $\VVV_1,\dots,\VVV_p$ of dimension $2n+1-r$ and
continuous maps $\Psi_i:\VVV_i\to\CCC^r$ satisfying (2) of the
lemma and such that $\Psi_1(\VVV_1)\cup\dots\cup\Psi_p(\VVV_p)$ is
a neighborhood of $\bK$ in $\CCC^r$.

\noindent{\bf Step 2.} Let $A=[u,v]\subset\RR$ be a compact
interval such that $A\cap\inter h(Z')\neq\emptyset$ and let
$\bK\in\CCC^0(A)$. We claim that there exist connected open
manifolds $\VVV_1,\dots,\VVV_r$ of dimension equal to
$2n-1+\dim\JJ$ and continuous maps $\Psi_i:\VVV_i\to\CCC(A)$
satisfying: (1) the union
$\Psi_1(\VVV_1)\cup\dots\cup\Psi_r(\VVV_r)$ is a neighborhood of
$\bK$ in $\CCC(A)$, (2) for any $i$ the composition
$\pi_J\circ\Psi_i:\VVV_i\to\JJ$ is a smooth map, and (3) the map
$(b\circ\Psi_i,e\circ\Psi_i):\VVV_i\to h^{-1}(u)\times h^{-1}(v)$
is a smooth submersion for each $i$. Except from (3) everything
follows as in the proof of Lemma
\ref{lemma:no-critical-value-parametrizing}, and (3) is a
consequence of Lemma \ref{lemma:de-onto}. Similar statements hold
replacing $[u,v]$ by $(u,v]$ and $[u,v)$, replacing the map
$(b\circ\Psi_i,e\circ\Psi_i)$ by $b\circ\Psi_i$ in the first case
and by $e\circ\Psi_i$ in the second one, and decreasing the
dimensions of $\VVV_i$ one unit.

\noindent{\bf Step 3.} Let again $A=[u,v]\subset\RR$ be a compact
interval, and let $\bK=(J,(K,\dots))\in\CCC^1(A)$. We claim that
there exist connected open manifolds $\VVV_1,\dots,\VVV_r$ of
dimension equal to $2n-2+\dim\JJ$ and continuous maps
$\Psi_i:\VVV_i\to\CCC(A)$ satisfying: (1) the union
$\Psi_1(\VVV_1)\cup\dots\cup\Psi_r(\VVV_r)$ is a neighborhood of
$\bK$ in $\CCC(A)$, (2) for any $i$ the compositions
$\pi_J\circ\Psi_i:\VVV_i\to\JJ$ and
$(b\circ\Psi_i,e\circ\Psi_i):\VVV_i\to h^{-1}(u)\times h^{-1}(v)$
are smooth maps. This can be proved as in the proof of Lemma
\ref{lemma:no-critical-value-parametrizing} using the (un)stable
manifold theorem for critical sets on normally hyperbolic vector
fields (see e.g. Theorem 4.1 in \cite{HPS}). To be more concrete,
suppose for simplicity that $A\cap h(Z')=\emptyset$, let
$\bK=(J,(K,\dots))\in\CCC^1(A)$ and let $F'\subset F$ be the
connected component to which $K\cap F$ belongs. Let
$\pi_{\pm}:W_{\pm}\to F'$ be the (un)stable manifolds and the
corresponding submersions for the vector field $-I\XXX$. Then one
can identify a neighborhood of $\bK$ in $\CCC^1(A)$ with
$(W_+\times_{F'}W_-)\cap h^{-1}(u)\times h^{-1}(v).$

\noindent{\bf Step 4.} Now let $r\geq 1$ be any integer and let
$\bK=(J,(K,\dots))\in\CCC^r$. Assume that $h(b(\bK))<h(e(\bK))$,
the other cases (either the opposite inequality or equality) being
analogous. Let $K\cap F=\{f_1,\dots,f_r\}$, labelled in such a way
that $d_1=h(f_1)<\dots<d_r=h(f_r)$. Let $\eta>0$ be small enough
so that each $A_j:=[d_j-\eta,d_j+\eta]$ is disjoint from $h(Z')$.
Assume also for simplicity that $h(b(\bK))<d_1-\eta$ and
$h(e(\bK))>d_r+\eta$. Define the intervals
$A'_0=(h(b(\bK))-\eta,d_1-\eta]$, $A'_r=[d_r+\eta,h(e(\bK))+\eta)$
and, for any $1\leq j\leq r-1$, $A'_j=[d_{j-1}+\eta,d_j-\eta]$.
Let $X_{j,\pm}=h^{-1}(d_j\pm\eta)$. Then the following fiber
product gives a neighborhood of $\bK$ in $\CCC^r$:
$$\CCC^0(A_0')\times_{X_{1,-}}\CCC^1(A_1)\times_{X_{1,+}}
\CCC^0(A_1')\times_{X_{2,-}}\CCC^1(A_2)\times
\dots\times\CCC^1(A_r)\times_{X_{r,+}}\CCC^0(A_r').$$ Here the
fiber product is defined using the maps $(b,e):\CCC^1(A_j)\to
X_{j,-}\times X_{j,+}$ and $(b,e):\CCC^0(A_j)\to X_{j,+}\times
X_{j+1,-}$ for $1\leq j\leq r-1$, the cases $j=0,r$ being the
obvious generalizations. Using the results in Steps 2 and 3 and a
simple computation with dimensions, the result follows.
\end{pf}

By the same argument as in Step 1 of the proof of Lemma
\ref{lemma:broken-lines}, one can choose a countable set
$\{\bK_\nu\}\subset\CCC^0$ such that, denoting by
$\{\UUU_{\nu,q}\}_{q\in\QQQ_{\nu}}$ the manifolds and by
$\Phi_{\nu,q}:\UUU_{\nu,q}\to\CCC$ the maps constructed in Lemma
\ref{lemma:no-critical-value-parametrizing} for $\bK=\bK_{\nu}$,
the images $\Phi_{\nu,q}(\UUU_{\nu,q})$ cover $\CCC^0$. Let also
$\SSS_{\nu,q,l,i}\subset\UUU_{\nu,q}$ denote the submanifolds
given by the lemma.

Let $\MMM=\DD\times M$. Let $g:N\to X^2_{\Lambda}$ a smooth map,
where $N$ is a smooth manifold of dimension $\leq 2n-4$,
satisfying $\Omega_f\subset g(N)$. Define $\NNN=\DD\times N$.
Consider the maps
$$F:\MMM\to X^2_{\Lambda}\qquad\qquad
G:\NNN\to X^2_{\Lambda}$$ defined as $F(\xi,m)=\xi\circ f(m)$ and
$G(\xi,n)=\xi\circ g(n)$. Define also $\PPP=S^1\times\CCC$ and
$\PPP_{\Lambda}=E_{\Lambda}\times_{S^1\times S^1}\PPP$, and let
the map $\Theta_{\Lambda}:\PPP_{\Lambda}\to X^2_{\Lambda}$ be
defined generalizing in the obvious way the map
$\Theta_{J,\Lambda}$ in (\ref{eq:def-Theta-J}). Choose a covering
of $B_{\Lambda}$ by open sets $\{\ZZZ_{\lambda}\}$ in such a way
that there exist trivializations
$E_{\Lambda}|_{\ZZZ_{\lambda}}\simeq\ZZZ_{\lambda}\times S^1$, and
denote by $\zeta_{\lambda}:\ZZZ_{\lambda}\times
S^1\times\CCC\to\PPP_{\Lambda}|_{\ZZZ_{\lambda}}$ the induced
trivializations of the bundle $\PPP_{\Lambda}$. Both maps $F$ and
$G$ submersions, so they are transverse to the following
compositions of smooth maps:
$$\xymatrix{e_{{\lambda},\nu,q}:\ZZZ_{\lambda}\times S^1\times\UUU_{\nu,q}
\ar[rrr]^-{\Id\times\Id\times\Phi_{\nu,q}} &&&
\ZZZ_{\lambda}\times S^1\times\CCC\ar[r]^-{\zeta_{\lambda}}
&\PPP_{\Lambda}\ar[r]^-{\Phi_{\Lambda}}&X^2_{\Lambda},}$$ and
$$\xymatrix{e'_{{\lambda},r,i}:\ZZZ_{\lambda}\times \{1\}\times\VVV_{r,i}
\ar[rrr]^-{\Id\times\iota\times\Psi_{r,i}} &&&
\ZZZ_{\lambda}\times S^1\times\CCC\ar[r]^-{\zeta_{\lambda}}
&\PPP_{\Lambda}\ar[r]^-{\Phi_{\Lambda}}&X^2_{\Lambda},}$$ where
$\iota:\{1\}\to S^1$ is the inclusion. For the same reason $F$ and
$G$ are transverse to the restriction of $e_{l,\nu,q}$ to each of
the manifolds $\ZZZ_{\lambda}\times S^1\times \SSS_{\nu,q,l,i}$.
Hence we have six countable sequences of smooth manifolds:
$\CS(e_{\lambda,\nu,q},F)$,
$\CS(e_{\lambda,\nu,q}|_{\ZZZ_{\lambda}\times
S^1\times\SSS_{\nu,q,l,i}},F)$, $\CS(e'_{\lambda,r,i},F)$, and the
same ones replacing $F$ by $G$. We denote by $\bCS$ the collection
of all these manifolds. Each of the manifolds in $\bCS$ projects
smoothly to $\JJ\times\DD$, and by Sard's theorem there exists a
residual subset $\Omega\subset\JJ\times\DD$ of regular values of
all these maps. Furthermore,
\begin{align*} \dim
\CS(e_{\lambda,\nu,q},F) &= \dim\JJ+\dim\DD, \\
\dim \CS(e_{\lambda,\nu,q}|_{\ZZZ_{\lambda}\times
S^1\times\SSS_{\nu,q,l,i}},F) &= \dim\JJ+\dim\DD-1\qquad\text{for
$i=1,2$},
\end{align*}
and all the remaining manifolds in $\bCS$ have dimension
$\leq\dim\JJ+\dim\DD-2$.

\begin{lemma}
\label{lemma:clau} Define for any $J\in\JJ$ the preimage
$\VVV_{J,r,i}:=(\pi_J\circ \Psi_{r,i})^{-1}(J)$. The omega-limit
set $\Omega_{\Theta_{J,\Lambda}}$ is contained in the union of the
sets $e_{\lambda,r,i}'(\ZZZ_i\times\{1\}\times \VVV_{J,r,i})$.
\end{lemma}
\begin{pf}
Let $\{\bx_i\}\subset\PPP_{J,\Lambda}$ be a diverging sequence
such that $\Theta_{J,\Lambda}(\bx_i)$ converges in
$X^2_{\Lambda}$. Recall that $p:X^2_{\Lambda}\to B_{\Lambda}$
denotes the projection. Passing to a subsequence we may assume
that $\{p(\bx_i)\}\subset\ZZZ_{\lambda}$ for some $\lambda$, so
that we can write
$\zeta_{\lambda}^{-1}(\bx_i)=(\beta_i,\theta_i,\bK_i)\in\ZZZ_{\lambda}\times
S^1\times\CCC_J$. Passing again to a subsequence we may assume
that $\beta_i\to \beta$, $\theta_i\to\theta$ and $\bK_i\to\bK$.
Suppose that $\bK=(J,(K,K_1,\dots,K_s,b))$. Since $\{\bx_i\}$
diverges, $\bK$ does not belong to $\CCC_J^0$ and consequently
$\{y_1,\dots,y_r\}:=K\cap F\neq\emptyset$, where we may suppose
that $h(y_1)<\dots< h(y_r)$. There are two cases to consider,
either $h(b)\leq h(y_1)$ or $h(y_r)\leq h(b)$. In the first case
consider the (discontinuous) map $\rho:X\to X$ defined by
$\rho(x)=x$ if $h(x)>h(y_1)$ and $\rho(x)=\theta\cdot x$ if
$h(x)\leq h(y_1)$. This map lifts to maps $\rho:R_i^{\sharp}\to
R_i^{\sharp}$. Define $\bK'=(J,(K',K'_1,\dots,K'_s,b'))$ by
setting $K'=\rho(K)$, $K_i'=\rho(K_i)$, $b'=\rho(b)=\theta\cdot
b$. It turns out that $\bK'\in\CCC_J$ and that
$\Phi_{J,\Lambda}\circ\zeta_l(\beta,\theta,\bK)=
\Phi_{J,\Lambda}\circ\zeta_l(\beta,1,\bK')$, so the result follows
from (1) in Lemma \ref{lemma:broken-lines}. The case $h(y_r)\leq
h(b)$ is dealt with similarly.
\end{pf}

Combining the previous lemma with the estimates above on the
dimensions of the manifolds in $\bCS$, together with Lemma
\ref{lemma:CS-compacte}, it follows that for $(J,\xi)\in\RRR$ the
space $\TTT_{J,\xi}$ is a zero dimensional compact manifold. This
proves (1) of Theorem \ref{thm:definition-pairing}. Claims (2) and
(3) also follow from the estimate on the dimension. To prove Claim
(4), suppose that $(J,\xi),(J',\xi')\in\RRR$. By standard
arguments, there exists a smooth path
$\gamma:[0,1]\to\JJ\times\DD$ going from $(J,\xi)$ to $(J',\xi')$
which is transverse to the projections to $\JJ\times\DD$ from each
of the manifolds in $\bCS$. This implies (using again the
estimates on dimensions, Lemma \ref{lemma:CS-compacte} and Lemma
\ref{lemma:clau}) that
$$\TTT=\{(x,y,t)\in \PPP_{\Lambda}\times\MMM\times [0,1]\mid
\pi(x,y)=\gamma(t), \Theta_{\Lambda}(x)=F(y) \}$$ is compact
oriented graph (here $\pi:\PPP_{\Lambda}\times\MMM\to\JJ\times\DD$
is the projection). More precisely, there is a decomposition
$\TTT=\TTT_{\edge}\cup\TTT_{\vertex}\cup\TTT_{\partial},$ where
$$\TTT_{\vertex}=\{(x,y,t)\in\TTT\mid
x\in\zeta_{\lambda}(\ZZZ_{\lambda}\times
S^1\times\SSS_{\nu,q,l,i})\text{ for some $\lambda,q,l$ and
$i\in\{1,2\}$ }\},$$ and hence, by transversality, is a finite
set, and
$$\TTT_{\partial}=\{(x,y,t)\in\TTT\mid t\in\{0,1\}\}=
\TTT_{J,\xi}\cup \TTT_{J',\xi'}.$$ Finally,
$\TTT_{\edge}=\TTT\setminus(\TTT_{\vertex}\cup\TTT_{\partial})$ is
an oriented 1-manifold, so each of its connected components (which
we call edges) $\gamma$ has a beginning and an end,
$\textbegin(\gamma),\textend(\gamma)\in\TTT_{\vertex}\cup\TTT_{\partial}$.
For each $p\in\TTT_{\vertex}\cup\TTT_{\partial}$ there is a
positive integer $k$, which is equal to $1$ if and only if
$p\in\TTT_{\partial}$, and a neighborhood of $p$ in $\TTT$
homeomorphic to a neighborhood of $0$ in $\{z\in\CC\mid
z^k\in\RR_{\geq 0}\}$. In particular any $p\in\TTT_{\partial}$ is
an extreme (either beginning or end) of a unique edge. Denote, for
any $p=(x,y,t)\in\TTT_{\partial}$, $\sigma(p)=\sigma(x,y)$.
Reversing the orientation of $\TTT$ if necessary, we may assume
the following. For any $p\in\TTT_{J,\xi}\subset\TTT_{\partial}$
belonging to the extremes of $\gamma$, $\sigma(p)=1$ if
$p=\textbegin(\gamma)$, and $\sigma(p)=-1$ if
$p=\textend(\gamma)$; similarly, for any
$p'\in\TTT_{J',\xi'}\subset\TTT_{\partial}$ belonging to the
extremes of $\gamma'$, $\sigma(p')=-1$ if
$p'=\textbegin(\gamma')$, and $\sigma(p')=1$ if
$p'=\textend(\gamma')$. Denote, for any
$p=(x,y,t)\in\TTT_{\edge}$, $\weight(p)=\weight(x)$. This defines
a locally constant function on $\TTT_{\edge}$ so for any
$\gamma\in\pi_0(\TTT_{edge})$ we have $\weight(\gamma)\in\QQ$. For
any $p\in\TTT_{\vertex}$ we have
$\sum_{\textbegin(\gamma)=p}\weight(\gamma)=
\sum_{\textend(\gamma)=p}\weight(\gamma),$ where both sums run
over the set of edges. Putting together all these observations, we
deduce that
$$\sum_{(x,y)\in\TTT_{J,\xi}}
\sigma(x,y)\weight(x)=\sum_{(x,y)\in\TTT_{J',\xi'}}
\sigma(x,y)\weight(x),$$ which is what we wanted to prove. The
same ideas allow to prove that $\Delta_{\Lambda}(\beta)$ is
independent of the chosen pseudocycle $f:M\to X^2_{\Lambda}$.
Finally, if two subspaces $D,D'\subset
C^{\infty}(X^2_{\Lambda},TX^2_{\Lambda})$ satisfy the requirements
of the theorem, then so does $D+D'$, and this allows to prove that
the definition of $\Delta_{\Lambda}(\beta)$ is independent of the
choice of $D$.

\section{Proofs of Theorem \ref{thm:classe-diagonal},
Corollary \ref{cor:existence-inverse} and Theorem
\ref{thm:singular-m}} \label{s:proof-thm-S1}

\subsection{Proof of Theorem \ref{thm:classe-diagonal}}
\label{ss:proof:thm:classe-diagonal} Without loss of generality we
can assume that $m$ is not contained in $h(Z')$. Let $\Lambda$ be
big enough so that for any $\Lambda'\geq\Lambda$ the inclusion
$\iota_{\Lambda'}:X_{\Lambda'}^2\to X_{\Lambda'+1}^2$ induces an
isomorphism between $(2n-2)$-dimensional cohomology groups. Then
we can identify $H^{2n-2}_{S^1\times S^1}(X\times X)$ with
$H^{2n-2}(X^2_{\Lambda})$ in such a way that $[\Delta_{\CC^*}]$
corresponds to $[\Delta_{\CC^*}]_{\Lambda}$. Let $J\in\JJ$ be
arbitrary. Pick some small $\epsilon>0$ so that
$(m-\epsilon,m+\epsilon)$ is disjoint from $h(Z')$. Then the
subset
$$\PPP_{J,\Lambda}^{m,\epsilon}=\{x\in \PPP_J
\mid (h\times
h)(\Phi_{J,\Lambda}(x))\in(m-\epsilon,m+\epsilon)^2\}$$ carries a
natural structure of smooth manifold because any point
$x\in\PPP_{J,\Lambda}^{m,\epsilon}$ is contained in
$E_{\Lambda}\times_{S^1\times S^1}(S^1\times\CCC_J^{0,0})$ and (3)
in Theorem \ref{thm:atlas} (combined with local trivializations of
$E_{\Lambda}\to B_{\Lambda}$) provides local charts of
neighborhoods of $x$. Furthermore, $\PPP_{J,\Lambda}^{m,\epsilon}$
is an open neighborhood of $\PPP_{J,\Lambda}^{m}=H^{-1}(m,m)$,
where $H:\PPP_{J,\Lambda}^{m,\epsilon}\to
(m-\epsilon,m+\epsilon)^2$ is the map sending any $x\in
\PPP_{J,\Lambda}^{m,\epsilon}$ to $(h\times
h)(\Phi_{J,\Lambda}(x))$. Since $H$ is a submersion,
$\PPP_{J,\Lambda}^{m}$ is a smooth manifold and the map
\begin{equation}
\label{eq:equi-diag} \Phi_{J,\Lambda}:\PPP_{J,\Lambda}^m\to
h^{-1}(m)^2_{\Lambda}=E_{\Lambda}\times_{S^1\times
S^1}(h^{-1}(m)\times h^{-1}(m))
\end{equation} represents as a pseudocycle
the image of the cohomology class $[\Delta_{\CC^*}]_{\Lambda}$
under the restriction map $H^*(X^2_{\Lambda})\to
H^*(h^{-1}(m)^2_{\Lambda})$. The map (\ref{eq:equi-diag}) is an
immersion, because $h^{-1}(m)$ contains no fixed points (it fails
to be injective at the preimages of pairs $(x,y)\in
h^{-1}(m)^2_{\Lambda}$ where $x,y$ belong to an orbit whose
stabiliser is nontrivial). Hence the pseudocycle represented by
(\ref{eq:equi-diag}) can be identified with the Poincar\'e dual
(PD) of the homology class represented by
$$\Delta_{m,\Lambda}=\{(x,y)\in h^{-1}(m)^2_{\Lambda}
\mid S^1\cdot x=S^1\cdot y\}.$$ Consider the map
$$\pi:h^{-1}(m)^2_{\Lambda}\to B_{\Lambda}\times Y_m\times Y_m\to
Y_m\times Y_m,$$ where the first map is induced by the quotient
$h^{-1}(m)\to h^{-1}(m)/S^1=Y_m$ and the latter map is the
projection, and denote by
$$f:H^{2n+2}(Y_m\times Y_m)\to H^{2n+2}(h^{-1}(m)^2_{\Lambda})$$
the morphism induced by $\pi$, which for $\Lambda$ big enough is
an isomorphism. Since $\pi$ is a submersion of orbifolds and we
can identify $\Delta_{m,\Lambda}=\pi^{-1}(\Delta_m)$, we have
$$\PD([\Delta_{m,\Lambda}])=f \PD([\Delta_m]).$$
That this standard argument in differential topology works in the
context of orbifolds follows from the realization of Poincar\'e
duality in terms of differentiable forms, as in \S 6 of \cite{BT}
(recall that the de Rham complex for orbifolds is defined in terms
of smooth invariant forms on local uniformizers which patch
together in the obvious sense). This finishes the proof of Theorem
\ref{thm:classe-diagonal}.

\subsection{Proof of Corollary \ref{cor:existence-inverse}}
\label{ss:cor:existence-inverse}

We will need the following result.
\begin{lemma}
\label{lemma:decomp}
 Take any decomposition $[\Delta_m]=\sum
e_i\otimes f^i\in H^*(Y_m)\otimes H^*(Y_m)$. For any cohomology
class $a\in H^*(Y_m)$ we have
\begin{equation}
\label{eq:PD-delta} \sum \left(\int_{Y_m}\alpha\cup
e_i\right)f^i=a.
\end{equation}
\end{lemma}
\begin{pf}
This is well known in the case of smooth manifolds. The same proof
as given for example in p. 127 of \cite{BT} translates word by
word to the context of orbifolds via the use of local
uniformizers.
\end{pf}

We now prove Corollary \ref{cor:existence-inverse}. Take any
decomposition
$$[\Delta_{\CC^*}]=\sum \epsilon_i\otimes\eta^i\in
H_{S^1}^*(X)\otimes H_{S^1}^*(X)$$ and let $a\in H^*(Y_m)$ be any
cohomology class. The Kirwan map is compatible with K\"unneth in
the following sense: given any class $\delta\in H^*_{S^1\times
S^1}(X\times X)$, if we write $\delta=\sum\alpha_i\otimes\beta^i$
using the decomposition (\ref{eq:Kunneth}), then
$\kappa_m^2(\delta)=\sum\kappa_m(\alpha_i)\otimes\kappa_m(\beta^
i)$. In particular, Theorem \ref{thm:classe-diagonal} implies that
\begin{equation}
\label{eq:diag} \sum
\kappa_m(\epsilon_i)\otimes\kappa_m(\eta^i)=[\Delta_m].
\end{equation}
It follows from the definition of $l$ that
$$l_m(a)=\sum_i
\left(\int_{Y_m}a\cup\kappa_m(\epsilon_i)\right)\eta^i.$$ Applying
$\kappa_m$ to both sides, taking into account (\ref{eq:diag}) and
using \ref{lemma:decomp}, we compute:
$$\kappa_m l_m(a) = \sum_i
\left(\int_{Y_m}a\cup\kappa_m(\epsilon_i)\right)\kappa_m(\eta^i)=a.$$
Hence Corollary \ref{cor:existence-inverse} is proved.

\subsection{Proof of Theorem \ref{thm:singular-m}}

We first prove that ${\kappa_m'}^2[\Delta_{\CC^*}]=[\Delta_m']$,
where $[\Delta'_m]$ is the Poincar\'e dual of the diagonal class
in $Y_m'$. For that it suffices to check that the vector field
$\nabla h$ is transverse to ${h'}^{-1}(m)$ and to apply the same
arguments as in \S \ref{ss:proof:thm:classe-diagonal}. The
definition of $h'$ given in \cite{LeTo} depends on some choices
(which do not affect the map $\kappa_m$) and we will prove the
required transversality when $h'$ is a small enough perturbation
of $h$ (we might need a smaller perturbation than \cite{LeTo}).
Away from a neighborhood of the fixed point set in $h^{-1}(m)$ the
function $h'$ coincides with $h$, so its $m$-level set is
transverse to $\nabla h$. So it suffices to look at a neighborhood
of some fixed point component $Y\subset h^{-1}(m)$.

Let $N\to Y$ be the normal bundle of the inclusion $Y\subset X$,
with its induced complex hermitian structure. The action of $S^1$
induces a splitting in complex subbundles $N=V^+\oplus
V^-=(V_1^+\oplus\dots\oplus V_s^+)\oplus (V_1^-\oplus\dots\oplus
V_k^-)$ and $S^1$ acts on $V_i^{\pm}$ with weight
$\pm\lambda_i^{\pm}$ for some positive integers
$\lambda_1^+,\dots,\lambda_1^-,\dots$. There are neighborhoods
$U\subset X$ of $Y$ and $U_N\subset N$ of the zero section of $N$
and a diffeomorphism $f:U_N\to U$ such that for any
$(v^+,v^-)=(v_1^+,\dots,v_s^+,v_1^-,\dots,v_k^-) \in U_N$ we have
$h\circ f(v^+,v^-)=m+\|v^+\|^2-\|v^-\|^2$, where $\|v^+\|^2=\sum_i
\lambda_i^+\|v_i^+\|^2$ and $\|v^-\|^2=\sum_j
\lambda_j^-\|v_j^-\|^2$, and also
$$\nabla h\circ
f(v^+,v^-)=2(\lambda_1^+v_1^+,\dots,
\lambda_s^+v_s^+,-\lambda_1^-v_1^-,\dots,-\lambda_s^-v_k^-)+
O(\|v^+\|^2+\|v^-\|^2).$$ Assume that for some $\delta>0$ the set
$\{(v^+,v^-)\in N \mid \|v^+\|^2+\|v^-\|^2<3\delta\}$ is contained
in $U_N$. Let $\rho:\RR\to\RR$ a smooth function with $\rho'\leq
0$, $\rho(t)=1$ for $t<\delta$ and $\rho(t)=0$ for $t>2\delta$.
Take $\epsilon\in\RR\setminus\{0\}$ with
$|\epsilon|<\min\{\sup|\rho'(t)|,\delta\}$, and require $\epsilon$
to be positive if and only if $\rk V^+\leq \rk V^-$. Then the
restriction of $h'$ to $U$ is defined as $m+\phi\circ f^{-1}$,
where
$$\phi(v^+,v^-)=\|v^+\|^2-\|v^-\|^2+\epsilon\rho(\|v^+\|^2+\|v^-\|^2).$$
Hence we need to prove that if $(v^+,v^-)$ satisfies
$\|v^+\|^2+\|v^-\|^2<3\delta$ and $\phi(v^+,v^-)=0$, then
$d\phi(v^+,v^-)(\nabla h\circ f)>0$. Using $\phi(v^+,v^-)=0$ one
computes
\begin{multline*}
\frac{1}{2}d\phi(v^+,v^-)(\nabla h\circ
f)=\left(\sum_i(\lambda_i^+)^2\|v_i^+\|^2+\sum_j(\lambda_j^-)^2\|v_j^-\|^2
\right)+ \\ +\epsilon\rho'(v^+,v^-)
\left(\sum_i(\lambda_i^+)^2\|v_i^+\|^2-\sum_j(\lambda_j^-)^2\|v_j^-\|^2
\right)+O(\|v^+\|^3+\|v^-\|^3).
\end{multline*}
The first expressions in big parenthesis is not less than the
second one, and nonzero if $\phi(v^+,v^-)=0$, so if the $O$ term
was absent then we would have $d\phi(v^+,v^-)(\nabla h\circ f)>0$.
Picking $\delta$ small enough and $|\epsilon|\leq \sup|\rho'|/2$,
the $O$ term will be smaller than the first two terms, so the
gradient will still be $>0$. This finishes the proof that
${\kappa_m'}^2[\Delta_{\CC^*}]=[\Delta_m']$.

Arguing as in the proof of Corollary \ref{cor:existence-inverse}
we deduce from ${\kappa_m'}^2[\Delta_{\CC^*}]=[\Delta_m']$ that
the element
$(PD'\otimes\Id)\circ(\kappa_m'\otimes\Id)[\Delta_{\CC^*}]$ (where
$PD':H^*(Y_m')\to H^{2n-2-*}(Y_m')^*$ denotes the Poincar\'e
duality map) corresponds to a map $l_m':H^*(Y_m')\to H^*_{S^1}(X)$
which is a right inverse of $\kappa_m'$. By an argument in linear
algebra $(f_H^*\otimes\Id)\circ
(PD'\otimes\Id)\circ(\kappa_m'\otimes\Id)[\Delta_{\CC^*}]$
corresponds to $l_m:=l_m'\circ f_H:IH^*(Y_m')\to H^*_{S^1}(X)$,
which is a right inverse of $\kappa_m=f_H^{-1}\circ\kappa_m'$.
Recall that $PD:IH^*(Y_m)\to IH^{2n-2-*}(Y_m)^*$ denotes the
Poincar\'e duality map. Since $f_H$ preserves the intersection
pairing we have $f_H^*\circ PD'=PD\circ f_H^{-1}$, so
\begin{align*}
(f_H^*\otimes\Id)\circ
(PD'\otimes\Id)\circ(\kappa_m'\otimes\Id)[\Delta_{\CC^*}] & =
(PD\otimes\Id)\circ(f_H^{-1}\otimes\Id)\circ(\kappa_m'\otimes\Id)
[\Delta_{\CC^*}] \\ &=(PD\otimes\Id)\circ(\kappa_m\otimes\Id)
[\Delta_{\CC^*}].
\end{align*}
This proves the theorem.

\subsection{An example}
\label{ss:example}

Let $\gamma:\RR\to X$ be a gradient flow line of $h$ such that
$$\lim_{t\to-\infty}h(\gamma(t))=\sup h
\qquad\text{and}\qquad \lim_{t\to\infty}h(\gamma(t))=\inf h.$$
Then $E=S^1\cdot\gamma(\RR)\subset X$ is an $S^1$-invariant
$2$-dimensional sphere embedded in $X$, which defines an
equivariant cohomology class $\alpha\in H^{2n-2}_{S^1}(X)$ (for
example, taking any finite dimensional approximation
$X_\Lambda=S_\Lambda\times_{S^1}X$, we define $\alpha_\Lambda\in
H^{2n-2}(X_\Lambda)$ as the Poincar\'e dual of
$S_\Lambda\times_{S^1}E$; making then $\Lambda$ go to $\infty$,
the classes $\alpha_\Lambda$ define a unique class $\alpha$). The
class $\alpha$ is independent of the choice of $\gamma$.

One can prove that for any regular value $m\in\imag\RR$ we have
$l_m(PD[\pt])=\alpha.$ This implies in particular that if
$m'\in\imag\RR$ is another regular value then the map
$$\kappa_{m'}\circ l_m:H^{2n-2}(Y_m)\to H^{2n-2}(Y_{m'})$$
is an isomorphism of vector spaces. An immediate consequence is
that $l_m$ is not in general a morphism of rings. Indeed, if $l_m$
were a morphism of rings then $\kappa_{m'}\circ l_m$ would also be
a morphism of rings. But one can construct examples in which
$Y_{m}$ is the blow up of $Y_{m'}$ at a point and in this case,
denoting by $\epsilon\in H^2(Y_m)$ be the Poincar\'e dual of the
exceptional divisor, it can be checked that $\kappa_{m'}\circ
l_m(\epsilon)=0$. However, $\epsilon^{n-1}\neq 0$, so
$\kappa_{m'}\circ l_m(\epsilon^{n-1})\neq 0$.

\section{Actions of compact tori of arbitrary dimension}
\label{s:high-dim}

We now sketch how to generalize the previous constructions in
order to prove Theorem \ref{thm:classe-diagonal-2}. Fix a subgroup
$S^1\simeq T_0\subset T$ such that the $T_0$-fixed point set
coincides with the $T$-fixed point set and take a basis
$\bu=\{u_1,\dots,u_q\}$ of $\tlie$. For each $l$ let $\XXX_{l}$
denote the vector field generated by the infinitesimal action of
$u_l$. The hypothesis on $T_0$ implies that for any two connected
components $F',F''$ of the zero set of $\XXX_{u_l}$ satisfying
$f'=\la\mu(F'),u_l\ra <f''=\la\mu(F''),u_l\ra$ there is some
$f'<a<f''$ such that the set $X_a=\{x\in X\mid \la
\mu(x),u_l\ra=a\}$ does not contain any $T_0$-fixed point. Indeed,
by hypothesis a $T_0$-fixed point is a zero of $\XXX_l$, and the
set of values of the function $\la\mu(\cdot),u_l\ra$ evaluated at
zeroes of $\XXX_l$ is finite.

Since $T_0$-stabilizers of the points in the level sets $X_a$ are
all finite, we can construct $T_0$-invariant multivalued
perturbations of the equation $\gamma'=-I\XXX_l$ supported near
the sets of the form $X_a$, just as we did in the case of actions
of the circle. Thus we get a finite dimensional space of
perturbations $\JJ$ and, for any $1\leq l\leq q$ and any
$J\in\JJ$, a space $\CCC_{l,J}$ parameterizing oriented chains of
$J$-perturbed gradient segments of $\la\mu,u_l\ra$, as in \S
\ref{ss:space}. We can also define, generalizing \S
\ref{ss:P-space}, the spaces
$$\CCC_J=\{(\bK_1,\dots,\bK_q)\in
\CCC_{1,J}\times\dots\times\CCC_{q,J}\mid
b(\bK_{i+1})=e(\bK_i)\},$$ $\PPP_J=T\times \CCC_J$, the maps
$(b,e):\CCC_J\to X\times X$ sending $(\bK_1,\dots,\bK_q)$ to
$(b(\bK_1),e(\bK_q))$ and $\Theta_J:\PPP_J\to X\times X$ sending
$(\theta,\bK_1,\dots,\bK_q)$ to $(\theta\cdot b(\bK_1),e(\bK_q))$,
and the action of $T\times T$ on $\PPP_J$ defined as
$(\alpha,\beta)\cdot(\theta,\bK_1,\dots,\bK_q)=
(\alpha\beta^{-1}\theta,\beta\cdot\bK_1,\dots,\beta\cdot\bK_q)$.
Define also for any natural $r\geq 0$ the set $\CCC_J^r$ as the
union, over all tuples $r_1,\dots,r_q$ of nonnegative integers
adding $r$, of the sets $\CCC_J\cap
(\CCC_{1,J}^{r_1}\times\dots\times\CCC_{q,J}^{r_q})$. Finally, let
$\CCC_J^{0,0}=\CCC_J\cap(\CCC_{1,J}^{0,0}\times\dots\times\CCC_{q,J}^{0,0})$
and define the weight of $(\bK_1,\dots,\bK_q)\in\CCC_J^{0,0}$ to
be the product of weights $\weight(\bK_1)\dots\weight(\bK_q)$.

As in the case of $S^1$ there are finite dimensional
approximations of the universal bundle $ET\times ET\to BT\times
BT$ of the form $E_{\Lambda}\to B_{\Lambda}$ (which are the $q$-th
Cartesian product of the corresponding fibrations for $S^1$), for
any natural number $\Lambda$, and we can consider the fiberwise
product $\PPP_{J,\Lambda}=E_{\Lambda}\times_{T\times T}\PPP_J$ and
the map $\Phi_{J,\Lambda}:\PPP_{J,\Lambda}\to X^2_{\Lambda}$.

Theorems \ref{thm:atlas} and \ref{thm:definition-pairing}
generalize straightforwardly to the present situation with a few
modifications which we now explain. The parametrization of
neighborhoods of $\CCC_J$ given by Theorem \ref{thm:atlas} will be
given by manifolds of dimension $2n+q$. Similarly one should
modify the dimensions given in Lemmata
\ref{lemma:no-critical-value-parametrizing} and
\ref{lemma:broken-lines} by adding $q-1$ in each case. Finally,
Lemma \ref{lemma:clau} should be generalized as follows. Let
$T_1\subset T$ be a subgroup such that $T=T_0\times T_1$. Then the
omega limit set of the map $\Phi_{J,\Lambda}:\PPP_{J,\Lambda}\to
X^2_{\Lambda}$ can be covered by the images of smooth maps with
domains of the form $\ZZZ_i\times T_1\times\CCC_J^1$, where
$\ZZZ_i$ has the same meaning as in Lemma \ref{lemma:clau}. This
still implies that the omega limit set has codimension $\geq 2$
and hence allows to prove Theorem \ref{thm:definition-pairing}.

As in \S \ref{ss:def-Delta} one can make $\Lambda$ go to $\infty$
and obtain a cohomology class
$$[\Delta_\CC^{T_0,T}]\in H^{2n-2k}_{T\times T}(X\times X).$$
To check that this class is independent of the basis $\bu$ note:
(1) if we replace $\{u_1,\dots,u_k\}$ by $\{-u_1,\dots,u_k\}$ then
we get trivially the same cohomology class and (2) two different
basis are homotopic up to reversing orientation. By an easy
deformation argument similar to (4) in Theorem
\ref{thm:definition-pairing} we deduce the independence on $\bu$.

Finally, the same arguments as in the proof of Theorem
\ref{thm:classe-diagonal} (see \S \ref{s:proof-thm-S1}) allow to
prove that $\kappa^2_m([\Delta_{\CC^*}^{T_0,T}])=[\Delta_m]$

\end{document}